\documentclass[11pt]{amsart}
\usepackage{amsmath,amsthm,amsfonts,amscd,amssymb,eucal,latexsym}
\theoremstyle{definition}
\newtheorem{theorem}{Theorem}[section]
\newtheorem{corollary}[theorem]{Corollary}
\newtheorem{lemma}[theorem]{Lemma}    
\newtheorem{proposition}[theorem]{Proposition}
\newtheorem{remark}[theorem]{Remark}

\newtheorem{example}[theorem]{Example}

\newtheorem{definition}[theorem]{Definition}

\newcommand{\sa}{{\rm sa}}

\newcommand{\spn}{{\rm span}}
\newcommand{\rad}{{\rm rad}}
\newcommand{\re}{{\rm re}}

\newcommand{\diag}{{\rm diag}}
\newcommand{\id}{{\rm id}}

\newcommand{\coedist}{{\rm dist_{op}}}
\newcommand{\coedistp}{{\rm dist'_{op}}}
\newcommand{\Rcoedist}{{\rm dist}^R_{\rm op}}

\newcommand{\sdist}{{\rm dist_s}}
\newcommand{\qdist}{{\rm dist_q}}
\newcommand{\cqdist}{{\rm dist_{cq}}}
\newcommand{\Rcqdist}{{\rm dist}^R_{\rm cq}}
\newcommand{\dist}{{\rm dist}}
\newcommand{\cb}{{\rm cb}}
\newcommand{\Afn}{{\rm Afn}}
\newcommand{\ex}{{\rm ex}}
\newcommand{\fa}{{\rm fa}}
\newcommand{\UCP}{S}

\newcommand{\OM}{{\rm OM}}

\newcommand{\Cb}{{\mathbb C}}

\newcommand{\Rb}{{\mathbb R}}
\newcommand{\Nb}{{\mathbb N}}
\newcommand{\HOS}{{\rm HOS}}
\newcommand{\OS}{{\rm OS}}
\newcommand{\nudist}{{\rm dist_{nu}}}
\newcommand{\nudistp}{{\rm dist'_{nu}}}
\newcommand{\cB}{{\mathcal B}}
\newcommand{\cH}{{\mathcal H}}
\newcommand{\cK}{{\mathcal K}}
\newcommand{\cD}{{\mathcal D}}
\newcommand{\cE}{{\mathcal E}}
\newcommand{\rd}{{\rm rad}}
\newcommand{\cl}{{\rm c}}
\newcommand{\CAM}{C^*{\rm M}}
\newcommand{\CM}{{\rm CM}}
\newcommand{\Hdist}{{\rm dist_H}}
\newcommand{\GHdist}{{\rm dist_{GH}}}

\allowdisplaybreaks

\begin{document}

\title[Gromov-Hausdorff convergence]{On Gromov-Hausdorff convergence for
operator metric spaces}

\author{David Kerr}
\author{Hanfeng Li}
\address{\hskip-\parindent
David Kerr, Department of Mathematics, Texas A{\&}M University,
College Station TX 77843-3368, U.S.A.}
\email{kerr@math.tamu.edu}

\address{\hskip-\parindent
Hanfeng Li, Department of Mathematics, SUNY at Buffalo,
Buffalo NY 14260-2900, U.S.A.}
\email{hfli@math.buffalo.edu}

\date{February 12, 2007}

\begin{abstract}
We introduce an analogue for Lip-normed operator systems of the second
author's order-unit quantum Gromov-Hausdorff distance and prove that it
is equal to the first author's complete distance. This enables
us to consolidate the basic theory of what might be called
operator Gromov-Hausdorff convergence. In particular we establish a
completeness theorem and deduce continuity in quantum tori, Berezin-Toeplitz
quantizations, and $\theta$-deformations from work of the second author.
We show that approximability by Lip-normed matrix algebras is equivalent
to $1$-exactness of the underlying operator space and, by applying a result of
Junge and Pisier, that for $n\geq 7$ the set of isometry classes of
$n$-dimensional Lip-normed operator systems is nonseparable. We also
treat the question of generic complete order structure.
\end{abstract}

\maketitle

\section{Introduction}
In \cite{GHDQMS} Marc Rieffel introduced a notion of quantum Gromov-Hausdorff
distance for compact quantum metric spaces, which are 
order-unit spaces equipped with a kind of generalized Lipschitz seminorm 
called a Lip-norm. One of the principal motivations was to build
an analytic framework for explaining the kinds of convergence of spaces
in string theory that involve changes of topology (see \cite{GHDQMS} for a
discussion). In addition to Rieffel's analogues of the Gromov completeness
and compactness theorems, there have been developed various convergence and
continuity results which apply for instance to
$\theta$-deformations \cite{OUQGHD}, quantum tori \cite{GHDQMS,AQT},
and Berezin quantization \cite{MACS} (see also \cite{OUQGHD}).

Given the $C^*$-algebraic nature of the examples of primary interest
and the fact that unital $C^*$-algebras are not determined by their order-unit
structure, Rieffel posed the problem of how to develop a version of quantum
Gromov-Hausdorff distance which would incorporate algebraic or matricial
information so as to be able to fully distinguish the underlying
noncommutative topology. Two different methods for doing this have been
independently proposed by the present two authors. Working in the
setting of Lip-normed operator systems (or what might be more suggestively
dubbed ``compact operator metric spaces'' in accord with Rieffel's
terminology), the first author defined a matricial version of quantum
Gromov-Hausdorff distance called {\it complete distance} which
formally elaborates on Rieffel's definition so as to bring the
matrix state spaces, and hence the complete order structure, into the
picture \cite{MQGHD}. The
second author meanwhile devised a strategy for quantizing
Gromov-Hausdorff distance which operates entirely at the ``function'' level,
in the spirit of noncommutative geometry. This versatile approach was
implemented in both the order-unit and $C^*$-algebraic contexts under the
terminology {\it order-unit} (resp.\ {\it $C^*$-algebraic}) {\it quantum
Gromov-Hausdorff distance} \cite{OUQGHD,CAQGHD} and affords many technical
advantages, as the continuity results in \cite{OUQGHD,CAQGHD} illustrate.

The immediate aim of the present paper is to show that these two approaches
become naturally reconciled in the framework of operator systems, in parallel
with what happens in the order-unit case \cite{OUQGHD}. More precisely, we
introduce an operator system version of order-unit quantum Gromov-Hausdorff
distance called ``operator Gromov-Hausdorff distance'' and
prove that it coincides with complete distance. In fact the
methods for treating the order-unit situation can be transferred to the
matricial order framework, and so our main task here is to
supply the necessary operator-system-theoretic input, including
material on amalgamated sums.
As a consequence of the equivalence of the two perspectives we
can speak unambiguously of ``operator Gromov-Hausdorff
convergence'' and the ``operator Gromov-Hausdorff topology''.

Exploiting the viewpoint of operator Gromov-Hausdorff distance we establish 
a completeness theorem and infer continuity in quantum tori, Berezin-Toeplitz
quantizations, and $\theta$-deformations by comparison with the second 
author's distance $\nudist$ from \cite{CAQGHD}. Fr\'{e}d\'{e}ric 
Latr\'{e}moli\`{e}re's result on quantum Gromov-Hausdorff approximation of 
quantum tori by finite-dimensional ones \cite{AQT} is also observed to be 
valid at the operator level.

Furthermore, we show that certain problems particularly pertinent to our
operator version of quantum metric theory can be resolved at or reduced to
the purely ``topological'' level of the operator space structure. More
specifically, we prove that a Lip-normed operator system is a limit of
Lip-normed matrix algebras if and only if it is $1$-exact as an operator
space, and, by invoking a result of Marius Junge and Gilles Pisier, that 
for $n\geq 7$ the set $\OM_n$ of
(isometry classes of) $n$-dimensional Lip-normed operator systems is
nonseparable. We thereby obtain answers to some questions about complete
distance that were left open in \cite{MQGHD}.

The organization of the paper is as follows. In Section~\ref{S-amalgsum} we
define and discuss amalgamated sums of operator spaces and systems. Operator
Gromov-Hausdorff distance is introduced in Section~\ref{S-oGHd}, and
amalgamated sums of operator systems are used here for showing that we obtain
a metric on the set of isometry classes of Lip-normed operator systems.
We also show the Lipschitz equivalence of operator Gromov-Hausdorff distance
and complete distance in Section~\ref{S-oGHd}.
Section~\ref{S-complete} treats the completeness theorem. In
Section~\ref{S-matrix} we record the continuity results which follow from
comparison with $\nudist$ as defined for Lip-normed unital $C^*$-algebras. 
We also give here our characterization of matrix
approximability for the operator Gromov-Hausdorff topology in terms of
$1$-exactness, as well as a similar characterization for $\nudist$. In the
latter case, however, quasidiagonality must be
added to $1$-exactness, and so we deduce that the operator Gromov-Hausdorff
topology is strictly weaker than the $\nudist$ topology on the set $\CAM$ of
isometry classes of Lip-normed unital $C^*$-algebras. We furthermore show that
matrix approximability for $C^*$-algebraic quantum Gromov-Hausdorff distance
is equivalent to the $C^*$-algebra being an MF algebra, so that the 
$C^*$-algebraic quantum Gromov-Hausdorff topology on $\CAM$ is neither weaker 
nor stronger than the operator Gromov-Hausdorff topology. 
Section~\ref{S-nonseparable} establishes the nonseparability of $\OM_n$
for $n\geq 7$. What we in fact
demonstrate is that for each $n\geq 7$ there is a nonseparable set of
$n$-dimensional Lip-normed operator systems which are all
isometric to each other as compact quantum metric spaces.
Finally, in Section~\ref{S-generic} we describe the generic complete
order structure of a Lip-normed $1$-exact operator system under operator 
Gromov-Hausdorff distance. 
%as well as the generic $C^*$-algebraic structure of 
%a Lip-normed quasidiagonal exact unital $C^*$-algebra under $\nudist$.

Terminology and notation related to quantum and operator metric spaces is
established at the beginning of Section~\ref{S-oGHd}. For general references
on operator spaces and systems see \cite{ER,Pau,Pis}. The notation 
$d_\cb$ as appears in Sections~\ref{S-matrix} and \ref{S-nonseparable} refers 
to completely bounded Banach-Mazur distance. In this
paper we will not assume operator spaces and systems to be complete.
\medskip

\noindent{\it Acknowledgements.} D. Kerr was supported by the Alexander 
von Humboldt Foundation and JSPS. He is grateful to Joachim Cuntz at the
University of M\"{u}nster and Yasuyuki Kawahigashi at the
University of Tokyo for hosting his stays during the respective
2003--2004 and 2004--2005 academic years. Part of this work
was carried out while H. Li was at the University of Toronto. 
We thank Gilles Pisier for pointing out Proposition~\ref{P-algebraic} and the 
referee for explaining that the linear isomorphism in this result also holds
in the operator system case (Proposition~\ref{P-algebraicosy}). We also thank
Nate Brown for suggesting the connection to MF algebras as developed
in the last part of Section~\ref{S-matrix}.

\section{Amalgamated sums of operator spaces and systems}\label{S-amalgsum}

To show that the operator Gromov-Hausdorff distance defined
in Section~\ref{S-oGHd} yields a metric on the set of isometry classes of
Lip-normed operator systems,
we will need the notion of an amalgamated sum of operator systems. Since
the methods for dealing with amalgamated sums of operator spaces and
operator systems are the same, we will first discuss the former.

In analogy with the full amalgamated free product of $C^*$-algebras
\cite{Bla}, we define the amalgamated sum of operator
spaces via a universal property in the category of operator spaces
with complete contractions as morphisms:

\begin{definition}\label{D-amalgsumos}
Given operator spaces $X,Y$, and $V$ and completely isometric
linear maps $\varphi_X : V\to X$ and $\varphi_Y : V\to Y$, the
{\it amalgamated sum of $X$ and $Y$ over $V$} is an operator
space $E$ with complete contractions $\psi_X : X\to E$ and
$\psi_Y : Y\to E$ satisfying the following:
\begin{enumerate}
\item $\psi_X \circ \varphi_X = \psi_Y \circ \varphi_Y$,

%\item $\psi_X (X) + \psi_Y (Y)$ is dense in $E$, and

\item whenever $F$ is an operator space and $\pi_X : X\to F$ and
$\pi_Y : Y\to F$ are complete contractions satisfying $\pi_X \circ
\varphi_X = \pi_Y \circ \varphi_Y$, there is a unique complete contraction
$\pi : E\to F$ such that $\pi_X = \pi\circ \psi_X$ and
$\pi_Y =\pi\circ \psi_Y$.
\end{enumerate}
We denote $E$ by $X +_V Y$.
\end{definition}

\noindent Clearly $X +_V Y$ is unique up to complete isometry, if it exists.

The case $V=0$ is discussed in Section~2.6 of \cite{Pis}.

For any operator space $X$ we denote by $C^* (X)$ the universal $C^*$-algebra
associated to $X$, meaning that there is a fixed completely isometric embedding
$\psi_X : X\hookrightarrow C^*(X)$ such that $\psi_X(X)$ generates $C^*(X)$ as
a $C^*$-algebra and for any $C^*$-algebra
$A$ and complete contraction $\varphi : X\to A$ there exists a (unique)
$^*$-homomorphism $\Phi : C^* (X) \to A$ with $\varphi = \Phi\circ \psi_X$
\cite[Thm.\ 3.2]{Pestov} \cite[Thm.\ 8.14]{Pis}.
%(cf. the definition of the universal non-self-adjoint operator
%algebra $OA(A)$ associated to $A$ in \cite[page 112]{Pisier}).
%Clearly $\psi_X$ is a completely isometric embedding, and so we can regard
%$X$ as an operator subspace of $C^* (X)$.
Identifying $X$ with $\psi(X)$ we may regard $X$ as an operator
subspace of $C^* (X)$. If $Y$ is a subspace of $X$, then the
associated $^*$-homomorphism $C^* (Y)\to C^* (X)$ is injective by
Wittstock's extension theorem. We have the following
result, which is easy to see.

\begin{proposition}\label{P-universal}
For any completely isometric linear maps $\varphi_X:V\rightarrow X$ and
$\varphi_Y:V\rightarrow Y$, the sum $X+Y$ in the full amalgamated
free product $C^* (X) *_{C^*(V)} C^* (Y)$ is an amalgamated sum
$X+_V Y$.
%The amalgamated sum $X +_V Y$ is precisely $X+Y$ inside the full amalgamated
%free product $C^* (X) *_{C^*(V)} C^* (Y)$.
Furthermore, $C^* (X +_V Y) = C^* (X) *_{C^*(V)} C^*(Y)$.
\end{proposition}

\noindent Proposition~\ref{P-universal} generalizes Theorem~8.15 of
\cite{Pis}.
%It also yields alternate proofs for
%Proposition~\ref{P-existamalgsumos}
%and Proposition~\ref{P-embedamalgsumos} below.

Since the canonical $^*$-homomorphisms $C^*(X)\rightarrow C^* (X)
*_{C^*(V)} C^* (Y)$ and $C^*(Y)\rightarrow C^* (X) *_{C^*(V)} C^*
(Y)$ are both faithful \cite{Bla}, we get:

\begin{corollary} \label{P-embedamalgsumos}
For any completely isometric linear maps $\varphi_X:V\rightarrow
X$ and $\varphi_Y:V\rightarrow Y$, the canonical maps
$X\rightarrow X+_VY$ and $Y\rightarrow X+_VY$ are both complete
isometries.
\end{corollary}

The following was pointed out to the second author by Gilles Pisier.

\begin{proposition}[Pisier]\label{P-algebraic}
Denote by $X +^{\rm a}_V Y$ the algebraic amalgamated sum
$(X\oplus Y)/\{(\varphi_X (v), -\varphi_Y (v)) : v\in V\}$. If $V$
is closed in both $X$ and $Y$, then the natural map $\varphi : X
+^{\rm a}_V Y\to X +_V Y$ is a linear isomorphism. If both $X$ and
$Y$ are complete, then so is $X +_V Y$.
\end{proposition}

\begin{proof} Obviously $\varphi$ is surjective.
When $V=0$, $\varphi$ is injective and the norm on $X +_0 Y$ is
the $\ell_1$-norm \cite[Sect.\ 2.6]{Pis}. Supposing that $V$ is
closed in both $X$ and $Y$, we have that $U=\{\varphi_X
(v)-\varphi_Y (v) \in X +_0 Y : v\in V \}$ is a closed subspace of
$X +_0 Y$. Clearly $X +_V Y$ is the quotient space $(X +_0 Y)/U$.
Thus $\varphi$ is always injective.

When both $X$ and $Y$ are complete the amalgamated sum $X +_V Y$ does not
change if we replace $V$ by its completion. Thus we may assume that $V$ is
also complete. Then $X +_V Y = (X +_0 Y)/U$ is complete.
\end{proof}

%\noindent The case $V=0$ in the above proposition is discussed in
%Section~2.6 of \cite{Pis}.

We now pass to the operator system setting. In this case the morphisms are
u.c.p.\ (unital completely positive) maps. Recall that a
{\it complete order embedding} is a completely positive isometry
$\varphi$ from an operator system $X$ to an operator system $Y$ such that
$\varphi^{-1} : \varphi (X) \to X$ is completely positive. A completely
positive map from $X$ to $Y$ is a complete order embedding if and only if it
is completely isometric. A {\it complete order isomorphism} is a surjective
complete order embedding.

\begin{definition}\label{D-amalgsumosy}
Given operator systems $X, Y$, and $V$ with unital complete order
embeddings $\varphi_X : V \to X$ and $\varphi_Y : V \to Y$, the
{\it amalgamated sum of $X$ and $Y$ over $V$} is an operator
system $E$ with u.c.p.\ maps $\psi_X : X \to E$ and $\psi_Y : Y
\to E$ satisfying the following:
\begin{enumerate}
\item $\psi_X \circ \varphi_X = \psi_Y \circ \varphi_Y$,

%\item $\psi_X (X) + \psi_Y (Y)$ is dense in $E$, and

\item whenever $F$ is an operator system and $\pi_X : X \to F$ and
$\pi_Y : Y \to F$ are u.c.p.\ maps satisfying
$\pi_X \circ \varphi_X = \pi_Y \circ \varphi_Y$, there is a unique u.c.p.\ map
$\pi : E \to F$ such that $\pi_X = \pi\circ \psi_X$ and
$\pi_Y = \pi\circ \psi_Y$.
\end{enumerate}
We denote $E$ by $X +_V Y$.
\end{definition}

\noindent Clearly $X +_V Y$ is unique up to unital complete order isomorphism,
if it exists.

%\begin{proposition}\label{P-existamalgsumosy}
%$X +_V Y$ exists.
%\end{proposition}
%
%\begin{proof}
%We can argue here as in the alternate proof given for
%Proposition~\ref{P-existamalgsumos}, noting in this case that the algebraic
%amalgamated sum $X +^{\rm a}_V Y$ has a natural involution.
%\end{proof}
%
%\begin{proposition}\label{P-embedamalgsumosy}
%The maps $X \to X +_V Y$ and $X \to X +_V Y$ are unital complete order
%embeddings.
%\end{proposition}
%
%\noindent The proof is similar to that of
%Proposition~\ref{P-embedamalgsumos}.

For any operator system $X$ we denote by $C^*_{\rm h} (X)$ the
universal unital $C^*$-algebra associated to $X$, meaning that
there is a fixed unital complete order embedding $\psi_X :
X\hookrightarrow C^*_{\rm h} (X)$ such that $\psi_X (X)$ generates
$C^*_{\rm h} (X)$ as a $C^*$-algebra and for any unital $C^*$-algebra $A$
and u.c.p.\ map $\varphi : X\to A$ there exists a (unique)
$^*$-homomorphism $\Phi : C^*_{\rm h} (X) \to A$ with $\varphi = \Phi\circ
\psi_X$ \cite[Prop.\ 8]{KW}. Identifying $X$ with $\psi(X)$ we may
regard $X$ as an operator subsystem of $C^*_{\rm h} (X)$. If $Y$
is a subsystem of $X$, then the associated $^*$-homomorphism $C^*_{\rm h}
(Y)\to C^*_{\rm h} (X)$ is injective \cite[Prop.\ 9]{KW}. In parallel to
Proposition~\ref{P-universal} we have the following easily
verified facts.

%Giver an operator system $X$ we denote by $C^*_{\rm h} (X)$ its associated
%universal $C^*$-algebra. This means that there is a fixed complete contraction
%$\psi_X : X\to C^*_{\rm h} (X)$ such that for any $C^*$-algebra $A$ and
%u.c.p.\ map $\varphi : X \to A$ there exists a unique $^*$-homomorphism
%$\Phi : C^*_{\rm h} (X)\to A$ with $\varphi = \Phi\circ\psi_X$.
%The map $\psi_X$ is in fact a unital complete order embedding, and so we may
%view $X$ as an operator subsystem of $C^*_{\rm h} (X)$. If $Y$ is a subsystem
%of $X$, then the associated $^*$-homomorphism $C^*_{\rm h} (Y) \to
%C^*_{\rm h} (X)$ is injective by Arveson's extension theorem. In parallel to
%Proposition~\ref{P-universal} we have the following easily verified facts.

\begin{proposition}\label{P-universalosy}
For any unital complete order embeddings $\varphi_X:V\rightarrow
X$ and $\varphi_Y:V\rightarrow Y$, the sum $X+Y$ in the full
amalgamated free product $C^*_{\rm h} (X) *_{C^*_{\rm h}(V)}
C^*_{\rm h} (Y)$ is an amalgamated sum $X+_V Y$. Furthermore,
$C^*_{\rm h} (X +_V Y) = C^*_{\rm h} (X) *_{C^*_{\rm h}(V)}
C^*_{\rm h}(Y)$.
% The
%amalgamated sum $X +_V Y$ is precisely $X + Y$ inside the full
%amalgamated free product $C^*_{\rm h} (X) *_{C^*_{\rm h} (V)}
%C^*_{\rm h} (Y)$. Moreover, $C^*_{\rm h} (X +_V Y) = C^*_{\rm h}
%(X) *_{C^*_{\rm h} (V)} C^*_{\rm h} (Y)$.
\end{proposition}

\begin{corollary}\label{C-embedamalgsumosy}
For any unital complete order embeddings $\varphi_X:V\rightarrow
X$ and $\varphi_Y:V\rightarrow Y$, the canonical maps
$X\rightarrow X+_VY$ and $Y\rightarrow X+_VY$ are both unital
complete order embeddings.
\end{corollary}

%\noindent Applying Proposition~\ref{P-universalosy} we obtain alternate
%proofs of Propositions~\ref{P-existamalgsumos} and \ref{P-embedamalgsumosy}.

%We remark finally that we don't know whether the operator system
%analogue of Proposition~\ref{P-algebraic} holds.

Consider now the algebraic amalgamated sum $X +^{\rm a}_V Y$ defined as\linebreak
$(X\oplus Y)/\{(\varphi_X (v), -\varphi_Y (v)) : v\in V\}$. Suppose 
that $V$ is closed in both $X$ and $Y$. The referee has 
indicated to us that the linear isomorphism in Proposition~\ref{P-algebraic}
also holds in the operator system case by the following argument. 
A matrix ordered space structure may be put on $X +^{\rm a}_V Y$ by 
declaring that $(x+y)^* = x^* + y^*$ and that $z\in M_n (X +^{\rm a}_V Y )$ is
positive if for every $\varepsilon > 0$ there exist $x\in M_n (X)_+$ and
$y\in M_n (Y)_+$ such that $z + \varepsilon I_n = x + y$. The Choi-Effros axioms
for an abstract operator system \cite{CE} are readily verified, and we
thereby obtain an operator system satisfying the universal property in
Definition~\ref{D-amalgsumosy}. 
We thus have the analogue of Proposition~\ref{P-algebraic}:

\begin{proposition}\label{P-algebraicosy}
If $V$ is closed in both $X$ and $Y$, then the natural map $\varphi : X
+^{\rm a}_V Y\to X +_V Y$ is a linear isomorphism. If both $X$ and
$Y$ are complete, then so is $X +_V Y$.
\end{proposition}

We note finally that, as the following example demonstrates, the operator system 
and operator space amalgamated sums need not be canonically isometric despite our 
use of the common notation $X +_V Y$. 

\begin{example}\label{E-sum}
Let $X=Y=M_2 (\Cb )$. Let $V$ be the space of diagonal $2\times 2$ matrices
and let $\varphi_X:V\rightarrow X$ and $\varphi_Y:V\rightarrow Y$ be
the natural embeddings. Set 
\[ x =  \left[
\begin{matrix} 0  & 1 \\ 0 & 0 \end{matrix} \right] \in X, \quad y=\left[
\begin{matrix} 0  & 0 \\ 1 & 0 \end{matrix} \right]\in Y. \]
As indicated in the proof of Proposition~\ref{P-algebraic}, 
the norm of $(x, y)$ in the operator space amalgamated sum 
is equal to
\[ \inf_{v\in V} \big( \| x+v\|+\| y-v\| \big) = 2. \]
Set 
\[ p =  \left[
\begin{matrix} 1  & 0 \\ 0 & 0 \end{matrix} \right], \quad q=\left[
\begin{matrix} 0  & 0 \\ 0 & 1 \end{matrix} \right]. \]
Then 
\[
\left[ \begin{matrix} p  & x \\ x^* & q \end{matrix} \right] \in M_2(X)_+,
\quad \left[
\begin{matrix} q  & y \\ y^* & p \end{matrix} \right]\in M_2(Y)_+. \]
Thus 
\[
\left[
\begin{matrix} I  & (x,y) \\ (x,y)^* & I \end{matrix} \right]=
\left[
\begin{matrix} p  & x \\ x^* & q \end{matrix} \right]
+\left[
\begin{matrix} q  & y \\ y^* & p \end{matrix} \right]\in M_2(X+^{\rm a}_VY)_+. \]
Since $\left[
\begin{smallmatrix} I  & w \\ w^* & I \end{smallmatrix} \right]\in M_2(W)_+$ if and 
only if $\|w\|\le 1$ for an element $w$ in an operator system $W$ \cite[Prop.\ 1.3.2]{ER},
we see that the norm of $(x, y)$ in the operator system amalgamated sum  
is at most $1$. On the other hand, using the identity maps 
$X\rightarrow M_2(\Cb)$
and $Y\rightarrow M_2(\Cb)$ one sees 
that the norm of $(x, y)$ in the operator system amalgamated sum  
is at least $1$. 
Thus the norm of $(x, y)$ in the operator system amalgamated sum 
is exactly $1$.
\end{example}

\section{Operator Gromov-Hausdorff distance}\label{S-oGHd}

The proofs in this section and the next are modeled on ones from the
order-unit and $C^*$-algebraic cases in \cite{OUQGHD,CAQGHD}. In the
case that the argument directly translates and the aspects particular to the 
operator system setting have been dealt with, we will simply refer the
reader to \cite{OUQGHD} or \cite{CAQGHD}.

We begin by establishing some notation and terminology pertaining to
quantum and operator metric structures.

We write $\dist_{\rm H}^\rho$ and $\GHdist$ to designate
Hausdorff distance (with respect to a metric $\rho$, which will be omitted
when it is given by a norm on a linear space) and Gromov-Hausdorff distance,
respectively.

Given an operator system $X$ we denote by $\UCP_n (X)$ its $n$th matrix
state space, i.e., the set of all unital completely positive maps from $X$
into the matrix algebra $M_n$. We have a canonical identification
$\UCP_n (\overline{X} ) = \UCP_n (X)$ where $\overline{X}$ is the completion of $X$.
We write $X_\sa$ for the set of self-adjoint elements of $X$.

A {\it Lip-norm} on an order-unit space $A$ is a seminorm $L$ (which we will
allow to take the value $+\infty$) on $A$ such that
\begin{enumerate}
\item $L(a) = 0$ for all $a\in\Rb 1$, and

\item the metric $\rho_L$ defined on the state space $S(A)$ of $A$ by
$$ \rho_L (\sigma , \omega ) = \sup \{ | \sigma (a) - \omega (a) | : a\in A
\text{ and } L(a) \leq 1 \} $$
induces the weak$^*$ topology.
\end{enumerate}
Notice that $L$ must in fact vanish precisely on $\Rb 1$. 
The {\it closure} of $L$ is the Lip-norm $L^\cl$ on the completion
$\overline{A}$ given by
$$ L^\cl (a) = \inf \Big\{ \liminf_{n\to\infty} L(a_n ) : \{ a_n \}
\text{ is a sequence in } A \text{ with } \lim_{n\to\infty}
a_n = a \Big\} . $$
We say that $L$ is {\it closed} if $L = L^\cl$.
A {\it compact quantum metric space} \cite[Defn.\ 2.2]{GHDQMS} is a pair 
$(A,L)$ consisting of an order-unit space $A$ equipped with a Lip-norm $L$ (by
permitting infinite values we are not strictly observing the definition from
\cite{GHDQMS}, but this does not cause any problems since we can always
restrict to the subspace on which $L$ is finite if we wish).

Let $(A,L_A )$ and $(B,L_B )$ be compact quantum metric spaces.
We say that $(A,L_A )$ and $(B,L_B )$ are {\it isometric} if there is an
isometry from $(A,L_A )$ to $(B,L_B )$, i.e., a unital order
isomorphism $\varphi : \overline{A}\to \overline{B}$ such that 
$L^\cl_X = L^\cl_Y \circ \varphi$.
A Lip-norm $L$ on the order-unit direct sum $A\oplus B$ is said to
be {\it admissible} if it induces $L_A$ and $L_B$ under the natural quotient
maps onto the respective summands. The {\it quantum Gromov-Hausdorff distance}
between $(A,L_A )$ and $(B,L_B )$ is defined by
$$ \qdist (A,B) = \inf \dist_{\rm H}^{\rho_L} (S(A),S(B)) $$
where the infimum is taken over all admissible Lip-norms $L$ on $A\oplus B$.
This yields a metric on the set of isometry classes of compact quantum metric
spaces \cite[Thm.\ 7.8]{GHDQMS}.

By a {\it Lip-normed operator system} we mean a pair $(X,L)$ where
$X$ is an operator system and $L$ is a Lip-norm on $X_\sa$ (which is the same
as saying that $(X_\sa , L)$ forms a compact quantum metric space). We will
also speak of Lip-normed unital $C^*$-algebras, Lip-normed exact operator 
systems, etc., when we need to qualify or specialize the class of operator 
systems under consideration. The {\it closure}
of $(X,L)$ is the Lip-normed operator system $(\overline{X} ,L^\cl )$ where
$L^\cl$ is the closure of $L$ on $\overline{X}_\sa$. We say that $(X,L)$ is
{\it closed} if it is equal to it closure.

%By a {\it compact $C^*$-algebraic
%metric space} we mean a pair $(A,L)$ where $A$ is a $C^*$-algebra and
%$L$ an adjoint-invariant seminorm on $A$ (possibly taking the value $+\infty$)
%such that $L$ vanishes precisely on $\Cb$ and the metric $\rho_L$
%defined on the state space $S(A)$ as in the order-unit case above gives rise
%to the weak$^*$ topology (equivalently, $(A, L |_{A_\sa} )$ is a Lip-normed 
%unital $C^*$-algebra) \cite[Defn.\ 2.2]{CAQGHD}. The reason for this
%additional specialized definition in the $C^*$-algebra case is that we want
%to consider convergence questions in the presence of the Leibniz rule or
%generalizations thereof (see Section~\ref{S-complete}).

Let $(X,L)$ be a Lip-normed operator system. We denote by $\rd (X)$ its
radius, i.e., the common value over $n\in\Nb$ of the radii of the metrics
$$ \rho_{L,n} (\varphi , \psi ) = \sup \{ \| \varphi (x) - \psi (x) \| :
x\in X \text{ and } L(x) \leq 1 \} $$
defined on the respective matrix state spaces $\UCP_n (X)$
(see \cite[Prop.\ 2.9]{MQGHD}). The closure $(\overline{X} ,L^\cl )$ satisfies
$\rho_{L^\cl ,n} = \rho_{L,n}$ for every $n$ (cf.\ \cite[Sect.\ 4]{MSS}). 
The Lip-norm unit ball $\{ x\in X_\sa : L(x) \leq 1 \}$ will be denoted by
$\cE (X)$, or $\cE (X,L)$ in case of confusion. For $R\geq 0$ we set
$$ \cD_R (X) = \{ x\in X_\sa : L(x) \leq 1 \text{ and } \| x \| \leq
R \} , $$
and in the case $R = \rd (X)$ we will simply write $\cD (X)$. We will 
frequently use the fact that $\cE (X) = \cD (X) + \Rb 1$ 
\cite[Lemma 4.1]{OUQGHD}.

Let $(X, L_X )$ and $(Y, L_Y )$ be Lip-normed operator systems. The
{\it complete distance} between $(X, L_X )$ and $(Y, L_Y )$ is defined by
$$ \sdist (X,Y) = \inf \sup_{n\in\Nb} \dist_{\rm H}^{\rho_{L,n}}
(\UCP_n (X),\UCP_n (Y)) $$
where the infimum is taken over all admissible Lip-norms $L$ on
$(X\oplus Y)_\sa$ \cite[Defn.\ 3.2]{MQGHD}.
By an {\it isometry} from $(X, L_X )$ to $(Y, L_Y )$ we mean a unital
complete order isomorphism $\varphi : \overline{X} \to \overline{Y}$ such that
$L^\cl_X = L^\cl_Y \circ \varphi$ on $\overline{X}_\sa$. When there exists an
isometry from $(X, L_X )$ to $(Y, L_Y )$, we say that $(X, L_X )$ and
$(Y, L_Y )$ are {\it isometric}. We denote by $\OM$ the set of isometry
classes of Lip-normed operator systems, and by $\OM^R$ the
subset consisting of isometry classes of Lip-normed operator
systems with radius no bigger than $R$.
The set of isometry classes of Lip-normed unital $C^*$-algebras will be denoted
$\CAM$. Complete distance defines a metric on $\OM$ \cite[Thm.\ 4.10]{MQGHD}.

In analogy to Definition~4.2 of \cite{OUQGHD} and Definition~3.3 of
\cite{CAQGHD} we introduce the following notion of distance for
Lip-normed operator systems.

\begin{definition}\label{D-oGHd}
Let $(X, L_X )$ and $(Y, L_Y )$ be Lip-normed operator systems.
We define the {\it operator Gromov-Hausdorff distance}
$$ \coedist (X,Y) = \inf \dist_{\rm H} (h_X (\cE (X)) ,
h_Y (\cE (Y))) $$
where the infimum is taken over all triples $(V, h_X , h_Y )$ consisting of
an operator system $V$ and unital complete order embeddings $h_X : X\to V$ and
$h_Y : Y\to V$.
We also define 
$$ \coedistp (X,Y) = \inf \dist_{\rm H} (h_X (\cD (X)) ,
h_Y (\cD (Y))) $$
and, for $R\geq 0$, 
$$ \Rcoedist (X,Y) = \inf \dist_{\rm H} (h_X (\cD_R (X)) ,
h_Y (\cD_R (Y))) $$
with the infima being taken over the same set of triples.
\end{definition}

The distance $\coedistp$ is the more immediate analogue of 
order-unit and $C^*$-algebraic quantum Gromov-Hausdorff distance 
\cite{OUQGHD,CAQGHD}. However, since we are considering unital embeddings in 
our present context, we can remove the norm restriction to obtain the simpler
definition $\coedist$. Indeed we will show that $\coedist$ and $\coedistp$ 
define Lipschitz equivalent metrics on $\OM$.
The reason for the $R$ version is to facilitate the proof of completeness
(see Section~\ref{S-complete}) as well as some arguments
involving continuity in continuous fields (see \cite[Sect.\ 7]{OUQGHD} and
the beginning of Section~\ref{S-matrix} below).

Although $\cE (X)$ is not itself totally bounded, it is, as pointed out above, 
equal to the set of scalar translations of the totally bounded set 
$\cD (X)$, so that for any triple $(V, h_X , h_Y )$ 
as in Definition~\ref{D-oGHd} we have
$$ \dist_{\rm H} (h_X (\cE (X)) , h_Y (\cE (Y))) \leq
\dist_{\rm H} (h_X (\cD (X)) , h_Y (\cD (Y))) $$
and hence $\coedist \leq \coedistp$. It is also evident from the definitions
that for $R\geq \max (\rd (X) , \rd (Y) )$ we have 
$\coedist (X,Y)\leq \Rcoedist (X,Y)$.

To show that $\coedist$, $\coedistp$, and $\Rcoedist$ define metrics on $\OM$, 
we can proceed in the same manner as in Section~3 of \cite{CAQGHD}, granted 
that we have the appropriate operator
systems facts at hand. We will thus only explicitly indicate what we require
at the operator system level and refer the reader to Section~3 of \cite{CAQGHD}
for the main line of argument.

First we have the triangle inequality, 
%which can be established
%along the lines of the proofs of Lemma~4.6 in \cite{OUQGHD} and
%Lemma~3.6 in \cite{CAQGHD}.
for which we make use of amalgamated sums
of operator systems and in particular the fact that the summands
unitally complete order embed into the sum, as asserted by
Corollary~\ref{C-embedamalgsumosy}. More precisely, given any operator 
system $V_X$ (resp.\ $V_Z$) containing $X$ and $Y$ (resp.\ $Z$ and $Y$) as 
operator subsystems, we embed everything into the amalgamated sum
$V_X +_Y V_Z$ to obtain the desired estimate.

\begin{lemma}\label{L-triangle}
For any Lip-normed operator systems  $(X, L_X )$, $(Y, L_Y )$, and $(Z, L_Z )$
we have
\begin{align*}
\coedist (X,Z) &\leq \coedist (X,Y) + \coedist (Y,Z) , \\
\coedistp (X,Z) &\leq \coedistp (X,Y) + \coedistp (Y,Z) ,
\end{align*}
and, for $R\geq 0$,
$$ \Rcoedist (X,Z) \leq \Rcoedist (X,Y) + \Rcoedist (Y,Z) . $$
\end{lemma}

Secondly we must show that two Lip-normed operator systems are distance
zero apart if and only if they are isometric. For this we need to consider
direct limits of operator systems. Direct
limits of operator spaces are discussed on page 39 of \cite{ER}.
The following lemma is a direct consequence of the abstract
characterization of operator systems as matrix order unit spaces by Choi and
Effros \cite[Thm.\ 4.4]{CE}.

\begin{lemma}\label{L-directlimitosy}
Let  $\{ X_j\}_{j\in J}$
be an inductive system of operator systems where
the maps are unital complete order embeddings.
The algebraic inductive limit $\lim\limits_{\longrightarrow} X_j$
equipped with the natural $^*$-vector space structure, order unit,
matricial order structure, and matricial norms is an operator system.
\end{lemma}
%with unital complete order embeddings
%over an index set $J$
%Since we don't know whether or not operator systems are preserved under
%taking quotients we can't apply the proof given there to the operator system
%case.
%
%\begin{lemma}\label{L-directlimitosy}
%Let $X_1 \stackrel{\varphi_1}{\longrightarrow}
%X_2 \stackrel{\varphi_2}{\longrightarrow} \cdots$ be a directed system
%of operator systems where the maps are unital complete order embeddings. Let
%$\lim_{\to} X_j$ be the algebraic direct limit, regarded as a $^*$-vector
%space with a distinguished unit. Then
%$\lim_{\to} X_j$ equipped with the natural matrix norms is an operator system.
%\end{lemma}
%
%\begin{proof}
%Using Arveson's extension
%theorem we can find a Hilbert space $\cH$ and unital complete order embeddings
%$\psi_j : X \to \cB (\cH )$ for each $j\in \Nb$ such that
%$\psi_j = \psi_{j+1} \circ \varphi_j$ for all $j$. Then
%$\lim_{\to} X_j$ becomes a $^*$-linear subspace of $\cB (\cH )$. It is evident
%that the limit matrix norms on $\lim_{\to} X_j$ are precisely those induced
%from those of $\cB (\cH )$.
%\end{proof}
%
%Notice that above proof also works for the case of operator spaces.
%In general, given an inductive system $\{ X_j\}_{j\in J}$
%of operator systems with unital complete order embeddings
%over an index set $J$, the limit $\lim_{\to} X_j$
%is an operator system equipped with the natural order unit and
%order structure (and norms). This follows directly from the abstract
%characterization of operator systems as matrix order unit spaces by Choi and
%Effros \cite[Thm.\ 4.4]{CE}.

With Corollary~\ref{C-embedamalgsumosy} and
Lemma~\ref{L-directlimitosy} at our disposal we now can argue as
in the proof of Theorem~3.15 in \cite{CAQGHD} to deduce that
distance zero is equivalent to being isometric. In view of
Lemma~\ref{L-triangle} we thereby conclude the following.

\begin{theorem}\label{T-metricOM}
The distances $\coedist$, $\coedistp$, and $\Rcoedist$ define metrics on $\OM$
and $\OM^R$, respectively.
\end{theorem}

\begin{corollary}\label{C-metricCAM}
The restrictions of $\coedist$, $\coedistp$, and $\Rcoedist$ define metrics 
on $\CAM$ and $\CAM^R$, respectively.
\end{corollary}

Next we establish several inequalities involving $\coedist$, $\coedist'$, and
$\Rcoedist$, including comparisons with complete distance $\sdist$ as
introduced in \cite{MQGHD}.
Notice that for any operator system $X$ the
pairing between $X$ and $\UCP_n (X)$ gives us a natural u.c.p.\
map $X\to C(\UCP_n (X), M_n )$. We need the following well-known fact.

\begin{lemma}\label{L-osyrep}
Let $X$ be an operator system, and consider the $C^*$-algebraic direct product
$\prod_{n=1}^{\infty} C(\UCP_n (X), M_n )$, defined as
$$ \big\{ (a_n)_{n\in\Nb} : a_n \in C(\UCP_n (X), M_n ) \text{ for all } n
\text{ and }\sup\nolimits_{n\in\Nb} \| a_n \| < \infty \big\} . $$
Then the natural linear map $\varphi : X\to
\prod_{n=1}^{\infty} C(\UCP_n (X), M_n )$ is a unital complete order embedding.
\end{lemma}

\begin{proof}
Clearly $\varphi$ is a u.c.p.\ map. Say $X\subseteq \cB (\cH )$.
Then $M_m (X)\subseteq \cB \big( \bigoplus^m_{j=1} \cH \big)$. Let $w\in
M_m (X)$ and $\varepsilon > 0$ be given. Then we can find a
finite-dimensional subspace $\cK$ of $\cH$ such that $\| qwq \| > \| w \|
- \varepsilon$, where $q = \diag (p, \dots , p)$ and $p$ is the orthogonal
projection onto $\cK$. Say $\dim (\cK ) = n$. Notice that the map
$\psi : X\to \cB (\cK ) = M_n$ sending $x$ to $pxp$ is in $\UCP_n (X)$, and
we have
\begin{align*}
\| (\id_{M_m}\otimes\varphi )(w) \| &\geq
\| ((\id_{M_m}\otimes\varphi )(w))(\psi) \|
= \| (\id_{M_m}\otimes\psi )(w) \| \\
&= \| qwq \| > \| w \| - \varepsilon .
\end{align*}
Thus $\| (\id_{M_m}\otimes\varphi )(w) \| = \| w\|$, and hence $\varphi$ is a
complete isometry.
\end{proof}

\begin{theorem}\label{T-ops}
For any Lip-normed operator systems $(X, L_X )$ and $(Y, L_Y )$ we have\linebreak
$\coedist (X,Y) = \sdist (X,Y)$.
\end{theorem}

\begin{proof}
Apply the same argument as in the proof of Proposition~4.7 in \cite{OUQGHD},
only this time using Lemma~\ref{L-osyrep} and substituting Arveson's extension
theorem for the Hahn-Banach theorem.
\end{proof}

Using Theorem~\ref{T-ops} we get a new proof of Theorem~4.10(ii) in
\cite{MQGHD}.

The following is the analogue of Proposition~4.8 in \cite{OUQGHD} and
Proposition~3.9 in \cite{CAQGHD}.

\begin{proposition}\label{P-comparison}
For any Lip-normed operator systems $(X, L_X )$ and $(Y, L_Y )$ we have
\begin{gather}
| \rd (X) - \rd (Y) |\leq \GHdist (\cD (X), \cD (Y))
\leq \coedistp (X,Y)\leq \rd (X) + \rd (Y), \label{E-1} \\
| \coedistp (X,Y) - \dist_{\rm op}^{\rd (X)} (X,Y) | \leq 
| \rd (X) - \rd (Y) | , \label{E-2} \\
\coedistp (X,Y) \leq 3\hspace*{0.4mm} \coedist (X,Y) . \label{E-3}
\end{gather}
For $R\geq 0$ we also have
\begin{gather}
\Rcoedist (X,Y)\leq 2\hspace*{0.4mm} \coedist (X,Y) . \label{E-4}
\end{gather}
\end{proposition}

\begin{proof}
The proofs of (\ref{E-1}) and (\ref{E-2}) are similar to those of (5) and (6)
of Proposition~4.8 in \cite{OUQGHD}. By the definition of $\sdist$ one has
$\qdist (X,Y)\leq \sdist (X,Y)$. The inequality (\ref{E-3}) then
follows from (\ref{E-2}), (\ref{E-4}), and the fact that
$| \rd (X) - \rd (Y) | \leq \qdist (X,Y)$. The proof of (\ref{E-4})
parallels those of (8) of Proposition~4.8 in
\cite{OUQGHD} and (6) of Proposition~3.9 in \cite{CAQGHD}.
\end{proof}

%The next result is the analogue of Proposition~4.9 in \cite{OUQGHD}.

%\begin{proposition}\label{P-comparison2}
%Let $(X, L_X )$ and $(Y, L_Y )$ be Lip-normed operator systems and let
%$R\geq \max (\rd (X) , \rd (Y) )$. Then we have
%\begin{gather}
%\sdist (X,Y)\leq \Rcoedist (X,Y), \label{E-5} \\
%\sdist (X,Y)\leq 2\hspace*{0.4mm} \coedistp (X,Y). \label{E-6}
%\end{gather}
%\end{proposition}

%\begin{proof}
%The inequality (\ref{E-6}) follows immediately from (\ref{E-5}), (\ref{E-2}),
%and (\ref{E-1}), while (\ref{E-5}) is a consequence of 
%Theorem~\ref{T-ops} since we clearly have $\Rcoedist \geq \coedist$. 
%\end{proof}

Putting together Theorem~\ref{T-ops}, Proposition~\ref{P-comparison}, and 
the observations in the second paragraph after
Definition~\ref{D-oGHd}, we have proved:

\begin{theorem}\label{T-equivalent}
The metrics $\sdist$, $\coedist$, and $\coedistp$ are Lipschitz equivalent on 
$\OM$, and on $\OM^R$ they are Lipschitz equivalent to $\Rcoedist$.
%More precisely,
%\begin{gather*}
%\frac13 \coedist \leq \sdist \leq 2\hspace*{0.4mm} \coedist , \\
%\frac12 \Rcoedist \leq \sdist \leq \Rcoedist .
%\end{gather*}
%Also, the metrics $\sdist$ and $\coedist$
%(resp.\ on $\OM^R$)
%(resp.\ $\sdist$ and $\Rcoedist$)
\end{theorem}

As a consequence of Theorem~\ref{T-equivalent} we can speak
of {\it operator Gromov-Hausdorff convergence} and the
{\it operator Gromov-Hausdorff topology} without any ambiguity. Throughout
the rest of the paper we will typically use the operator Gromov-Hausdorff
distance $\coedist$ in the formulation of convergence and completeness results
with the tacit understanding that these apply equally well to the complete
distance $\sdist$, as well as to $\coedistp$.

We also have the following two facts, which can be established along the lines
of the proofs of Theorem~3.16 and Proposition~3.17, respectively, in
\cite{CAQGHD}. We denote by $\CM$ the set of isometry classes of compact
metric spaces, and for a compact metric space $(X, \rho)$ we write $L_\rho$
for the associated Lipschitz seminorm on $C(X)$.

\begin{theorem}\label{T-CGH}
The map $(X, \rho)\mapsto (C(X),  L_{\rho} )$ is a homeomorphism from
$(\CM , \GHdist )$ onto a closed subspace of
$(\OM , \coedist )$.
\end{theorem}

\begin{proposition}\label{P-opGH}
Let $(X, \rho_X )$ and $(Y, \rho_Y )$ be compact metric spaces. For
any $R\geq 0$ we have
\begin{gather*}
\Rcoedist (C(X), C(Y))\leq \GHdist (X,Y). 
\end{gather*}
\end{proposition}

\section{Completeness}\label{S-complete}

%Along the same line of proof as for the completeness part of
%\cite[Thm.\ 4.4]{Li2}, one can show:
We establish in this section a completeness theorem for operator
Gromov-Hausdorff distance. 
%complete distance in view Theorem~\ref{T-equivalent}.

One way to obtain a Lip-normed unital $C^*$-algebra is to restrict the Lip-norm
of a $C^*$-algebraic compact quantum metric space \cite[Defn.\ 2.2]{CAQGHD}.
The latter type of Lip-norm is a complex scalar version of
an order-unit Lip-norm which is defined on the whole $C^*$-algebra (but
possibly taking the value $+\infty$) and required to be adjoint invariant,
vanish precisely on $\Cb 1$, and induce the weak$^*$ topology on the state
space via the associated metric (defined the same way as in the order-unit 
case). Such Lip-norms appear
naturally in many examples (e.g., quantum metrics arising from ergodic actions
of compact groups \cite{MSACG}), and in the $C^*$-algebraic case one may wish
to study convergence questions in the presence of the Leibniz rule or
generalizations thereof such as the $F$-Leibniz property, which we now
recall (see \cite[Sect.\ 5]{MQGHD} \cite[Sect.\ 4]{CAQGHD}).

Let $F : \Rb^4_+ \to \Rb_+$ be a continuous function which is nondescreasing
with respect to the partial order on $\Rb^4_+$ under which
$(x_1 , x_2 , x_3 , x_4 ) \leq (y_1 , y_2 , y_3 , y_4 )$ if and only
$x_j \leq y_j$ for each $j$. A $C^*$-algebraic compact quantum metric space
$(A,L)$ is said to satisfy the {\it $F$-Leibniz property} if
$$ L(ab) \leq F(L(a),L(b), \| a \| , \| b \| ) $$
for all $a,b\in A$. Note that taking
$F(x_1 , x_2 , x_3 , x_4 ) = x_1 x_4 + x_2 x_3$ yields the Leibniz rule.

We write $\CAM_F$ for the subset of $\CAM$ whose elements come from
$C^*$-algebraic compact quantum metric spaces satisfying the $F$-Leibniz
property.

\begin{theorem}\label{T-complete}
The metric space $(\OM , \coedist )$ is complete. Let $F : \Rb^4_+\to \Rb_+$
be a continuous nondecreasing function. Then $(\CAM_F, \coedist )$ is also
complete.
\end{theorem}

\begin{proof}
Let $\{ (X_n, L_n )\}_{n\in\Nb}$ be a Cauchy sequence in $(\OM ,
\coedist )$. We may assume each $L_n$ to be closed. By (\ref{E-1})
we have $R := 1 + \sup_{n\in\Nb} \rd_{X_n} < +\infty$. Thus $\{
(X_n, L_n )\}_{n\in\Nb}$ is also a Cauchy sequence in $(\OM^R ,
\Rcoedist )$ by Theorem~\ref{T-equivalent}. To show that $\{ (X_n,
L_n )\}_{n\in\Nb}$ converges, it suffices to show that a
subsequence converges under $\Rcoedist$. Thus, passing to a
subsequence, we may assume that $\Rcoedist (X_n , X_{n+1} ) <
2^{-n}$ for all $n$. By Corollary~\ref{C-embedamalgsumosy} and
Lemma~\ref{L-directlimitosy} we can find a complete operator
system $V$ containing all of the $X_n$ as operator subsystems
such that $\Hdist (\cD_R (X_n ), \cD_R (X_{n+1} )) < 2^{-n}$ for
all $n$. Since $V$ is a complete metric space, the set of
non-empty closed compact subsets of $V$ is complete with respect
to Hausdorff distance. Denote by $W$ the limit of the sequence $\{
\cD_R (X_n) \}_{n\in\Nb}$ with respect to Hausdorff distance.

Since each $\cD_R (X_n )$ is $\Rb$-balanced
(i.e., $\lambda x\in \cD_R (X_n )$ for all $x\in\cD_R (X_n )$ and
$\lambda\in\Rb$ with $| \lambda |\leq 1$)
and convex, has radius $R$, and contains $0_V$ and $R \cdot 1_V$,
we can clearly say the same about $W$. Thus the set
$\Rb_+ \cdot W = \{ \lambda w : \lambda\in\Rb_+, w\in W \}$ is a real linear
subspace of $V_\sa$ containing $1_V$. Denote it by $B$, and denote by $X$ 
the closure of $B+iB$. Then $X$ is an operator subsystem of $V$. Notice that
Lemmas~4.8 and 4.9 of \cite{CAQGHD} hold in our current context, that is,
$(X,L)$ is a Lip-normed operator system with $\rd (X) \leq R$ and
$W = \cD_R (B)$, where $L$ is defined by
\begin{gather}
L(x) := \inf \big\{ \limsup\nolimits_{n\to \infty} L_n (x_n ) : x_n \in
X_n \text{ for all } n \text{ and } \lim\nolimits_{n\to\infty} x_n = x
\big\} \label{E-limit}
\end{gather}
for all $x\in X$. Now we have $\Rcoedist (X_n , X)\leq \Hdist (\cD_R (X_n ), \cD_R(X))
%=\Hdist (\cD_R (X_n ), W)
\to 0$ as $n\to\infty$. This proves the first assertion.

Now assume further that $(X_n, L_n)$ lies in $\CAM_F$ for all $n$.
Let $Z$ be a countable dense subset of $W+iW$.
Passing to a subsequence, we may assume that for any $x,y \in Z$ there exist
$x_n , y_n \in X_n$ for each $n$ such that $x_n \to x$, $y_n \to y$, and
$\{ x_n y_n \}_{n\in\Nb}$ converges to an element in $B+iB$. This is
sufficient for Lemma~4.7 of \cite{CAQGHD} to hold, so that if
$x,y \in X$ and we have $x_n , y_n \in X_n$ for each $n$ such that $x_n \to x$
and $y_n \to y$ as $n\to \infty$ then
$\{ x_n y_n \}_{n\in\Nb}$ converges to an element in $X$ and the limit
depends only on $x$ and $y$. Then we can define a product on $X$ by setting
$x\cdot y$ to be the limit of $\{ x_n y_n \}_{n\in\Nb}$. An argument similar
to that after Lemma~4.7 in \cite{CAQGHD} shows that $X$ becomes a unital
$C^*$-algebra. Since each $(X_n , L_n)$ satisfies the $F$-Leibniz property,
so does $(X, L)$ in view of (\ref{E-limit}). This proves the second assertion.
\end{proof}

Note that the second assertion of Theorem~\ref{T-complete} is equivalent
to Theorem~5.3 of \cite{MQGHD} in view of Theorem~\ref{T-ops}.
%Notice that Theorem~\ref{T-complete} improves Theorem~5.3 of
%\cite{MQGHD}.

%We don't know how to prove the analogue of Lemma~4.6 in \cite{CAQGHD}
%for Lip-normed operator systems, and thus we don't know how
%to prove the analogue of the compactness part of Theorem~4.4 in \cite{CAQGHD}.
%Neither do we know whether or not continuous fields of
%operator systems can always be subtrivialized at any point
%(see Question~7.3 and Proposition~7.4 in \cite{OUQGHD}). Hence
%we don't know how to prove the analogue of (3)$\Rightarrow$(1)
%in Theorem~1.2 of \cite{OUQGHD}, though this also implies
%the analogue of the compactness part of Theorem~4.4 in \cite{CAQGHD}.

\section{Continuity, $\nudist$, $\cqdist$, and matrix 
approximability}\label{S-matrix}

In Remark~5.5 of \cite{CAQGHD} the second author introduced a distance
$\nudist$ for $C^*$-algebraic compact quantum metric spaces. We will apply
the unital version of the definition
to Lip-normed unital $C^*$-algebras, keeping the same notation. Thus for
Lip-normed unital $C^*$-algebras $(A,L_A )$ and $(B,L_B )$ we set
\begin{align*} 
\nudist (A,B) &= \inf \dist_{\rm H} (h_A (\cE (A)) , h_B (\cE (B))) ,\\
\nudistp (A,B) &= \inf \dist_{\rm H} (h_A (\cD (A)) , h_B (\cD (B))) ,
\end{align*}
where the infima are taken over all triples $(D, h_A , h_B )$ consisting of
a unital $C^*$-algebra $D$ and unital faithful $^*$-homomorphisms 
$h_A : A\to D$ and $h_B : B\to D$. The notation reflects
the fact that these distances behave well for unital nuclear
$C^*$-algebras in the context of continuity problems. Applying the same 
arguments as in the operator system case 
(see Section~\ref{S-oGHd}), it can be shown that $\nudist (A,B)$ and 
$\nudistp (A,B)$ define Lipschitz equivalent metrics on the set $\CAM$ of 
isometry classes of Lip-normed unital $C^*$-algebras
(more precisely, $\nudist \leq \nudistp \leq 3\hspace*{0.4mm}
\nudist$), and so it suffices for
our purposes to work with the simpler definition $\nudist$. 
It is easily seen that everything in \cite[Sect.\ 5]{CAQGHD} pertaining to 
$\nudist$ works in the unital situation as well.

It follows from the definitions that for any Lip-normed unital $C^*$-algebras
$(A , L_A )$ and $(B , L_B )$ we have $\nudist (A, B)\geq \coedist (A, B)$.
Thus Theorems~5.2 and 5.3 in \cite{CAQGHD} hold with $\cqdist$ and
$\Rcqdist$ replaced by $\coedist$ and $\Rcoedist$, respectively, when $T$
is a compact metric space and each fibre $A_t$ is
nuclear. Since $C^*$-algebras admitting ergodic actions of compact
groups are automatically nuclear, we see that Theorem~5.11 and
Corollary~5.12 in \cite{CAQGHD} hold with $\cqdist$ replaced by $\coedist$
when $T$ is a compact metric space.

In particular, this gives us continuity with respect to $\nudist$ and 
$\coedist$ in quantum tori, Berezin-Toeplitz 
quantizations, and $\theta$-deformations (see \cite[Sect.\ 5]{CAQGHD}).
Corollary~2.2.13 and Proposition~3.1.4 in \cite{AQT} also enable
us to conclude approximation of quantum tori by finite quantum
tori under $\nudist$ and $\coedist$ (see Theorem~1.0.1 in
\cite{AQT}).

We turn now to the problem of matrix approximability. This will in
particular enable us to distinguish the $\coedist$ and $\nudist$ topologies
on the set of isometry classes of Lip-normed unital $C^*$-algebras.

Recall that an operator space $X$ is said to be {\it $1$-exact} if for every
finite-dimensional subspace $E\subseteq X$ and $\lambda > 1$ there is an
isomorphism $\alpha$ from $E$ onto a subspace of a matrix algebra such that
$\| \alpha \|_\cb \| \alpha^{-1} \|_\cb \leq\lambda$ (i.e., if $X$ is exact
with exactness constant $1$). This is equivalent to requiring that for
every $C^*$-algebra $A$ and closed two-sided ideal $I\subseteq A$ the natural
complete contraction
$$ (A \otimes_{\rm min} X) / (I \otimes_{\rm min} X) \to (A/I)
\otimes_{\rm min} X $$
is isometric (see \cite[Sect.\ 14.4]{ER} or \cite[Sect.\ 17]{Pis}).
An operator system is
said to be $1$-exact if it is $1$-exact as an operator space. 

\begin{lemma}\label{L-coe}
Let $X$ be a $1$-exact operator system. Let $\mathcal{H}$ be a Hilbert
space and $\iota : X\to\mathcal{B} (\mathcal{H} )$ a unital complete order
embedding. Then there is a net
$$ X \stackrel{\varphi_\lambda}{\longrightarrow} M_{n_\lambda}
\stackrel{\psi_\lambda}{\longrightarrow} \mathcal{B} (\mathcal{H} ) $$
of unital completely positive maps through matrix algebras
such that $\psi_\lambda \circ\varphi_\lambda$ converges pointwise to $\iota$.
\end{lemma}

\begin{proof}
Since $X$ is $1$-exact a standard application of Wittstock's extension
theorem produces a net
$$ X \stackrel{\varphi_\lambda}{\longrightarrow} M_{n_\lambda}
\stackrel{\psi_\lambda}{\longrightarrow} \mathcal{B} (\mathcal{H} ) $$
of completely contractive maps through matrix algebras such that $\psi_\lambda
\circ\varphi_\lambda$ converges pointwise to $\iota$.
Applying the construction in the proof of
Proposition~3.6 in \cite{PS} we can then produce a net
$$ X \stackrel{\varphi'_\lambda}{\longrightarrow} M_{m_\lambda}
\stackrel{\psi'_\lambda}{\longrightarrow} \mathcal{B} (\mathcal{H} ) $$
of completely positive contractive maps through matrix algebras such that
$\psi'_\lambda \circ\varphi'_\lambda$ converges pointwise to $\iota$ and
$\lim_\lambda \| \varphi'_\lambda (1) - 1 \| = 0$. We may assume that for all
$\lambda$ we have $\| \varphi'_\lambda (1) - 1 \| < 1/2$ and
$\| \psi'_\lambda \circ\varphi'_\lambda (1) - 1 \| < 1/2$ so that both
$\varphi'_\lambda (1)$ and $\psi'_\lambda (1)$ are invertible, which permits
us to define the maps
\begin{align*}
\varphi''_\lambda (\cdot ) &= \varphi'_\lambda (1)^{-1/2} \varphi'_\lambda
(\cdot ) \varphi'_\lambda (1)^{-1/2} , \\
\psi''_\lambda (\cdot ) &= \psi'_\lambda (1)^{-1/2} \psi'_\lambda
(\cdot ) \psi'_\lambda (1)^{-1/2} .
\end{align*}
We thereby obtain a net
$$ X \stackrel{\varphi''_\lambda}{\longrightarrow} M_{m_\lambda}
\stackrel{\psi''_\lambda}{\longrightarrow} \mathcal{B} (\mathcal{H} ) $$
of unital completely positive maps such that
$\psi''_\lambda \circ\varphi''_\lambda$ converges pointwise to $\iota$, as
desired.
\end{proof}

\begin{lemma}\label{L-supsystem}
Let $X$ be a separable operator system and $(Y,L)$ a Lip-normed 
finite-dimensional operator subsystem of $X$. Let $\varepsilon > 0$. Then 
there is a Lip-norm $L'$ on $X$ with respect to which we have 
$\Hdist (\cE (X), \cE (Y)) \leq\varepsilon$.
\end{lemma}

\begin{proof}
Take any Lip-norm $L''$ on $X_\sa$. Since $Y_\sa$ is 
finite-dimensional we may assume by scaling $L''$ if necessary that 
$L'' |_{Y_\sa} \leq L$. The finite-dimensionality of $Y_\sa$ also guarantees
the existence of a bounded projection $P : X_\sa \to Y_\sa$ (see for example
\cite[Lemma 3.2.3]{BK}). Define the seminorm $L'$ on $X$ by
$$ L' (x) = \max \left( L(P(x)), L'' (x) , \varepsilon^{-1} \| x - P(x) \| 
\right) . $$
Since $L'$ vanishes precisely on $\Rb 1$, dominates $L''$, and is 
finite on a dense subspace of $X_\sa$, we deduce that $L'$ is a Lip-norm, so 
that $(X,L' )$ forms a Lip-normed operator system. 

Now if $x\in\cE (X)$ then $P(x) \in\cE (Y)$ and 
$\| x - P(x) \| \leq\varepsilon$. Since
$\cE (Y) \subseteq \cE (X)$ (because $L' |_{Y_\sa} = L$), we thus 
conclude that the Hausdorff distance between $\cE (X)$ and 
$\cE (Y)$ is at most $\varepsilon$, as desired.
\end{proof}

\begin{theorem}\label{T-exact}
For a Lip-normed operator system $(X,L)$ the following are equivalent:
\begin{enumerate}
\item $X$ is $1$-exact,
\item for every $\varepsilon > 0$ there is a Lip-normed operator system
$(Y,L' )$ such that $Y$ is an operator subsystem of a matrix algebra and
$\coedist (X,Y) \leq\varepsilon$,
\item for every $\varepsilon > 0$ there is a Lip-normed matrix
algebra $(M_n ,L' )$ such that\linebreak $\coedist (X, M_n ) \leq\varepsilon$.
\end{enumerate}
\end{theorem}

\begin{proof}
(i)$\Rightarrow$(ii). By Lemma~\ref{L-coe} there is a unital
complete order embedding $\iota : X\to\mathcal{B} (\mathcal{H} )$ and a net
%\begin{gather*}
%\xymatrix{
%X \ar[r]^-{\varphi_\lambda} & B_\lambda \ar[r]^-{\psi_\lambda} & \mathcal{B}
%(\mathcal{H} )}
%\end{gather*}
$$ X \stackrel{\varphi_\lambda}{\longrightarrow} M_{n_\lambda}
\stackrel{\psi_\lambda}{\longrightarrow} \mathcal{B} (\mathcal{H} ) $$
of unital completely positive maps through matrix algebras
such that $\psi_\lambda \circ\varphi_\lambda$ converges pointwise to $\iota$.
In view of the equality of $\coedist$ and $\sdist$
(Theorem~\ref{T-ops}) we can now proceed as in the proof of
Proposition~3.10 in \cite{MQGHD} to obtain (ii).

(ii)$\Rightarrow$(iii). Apply Lemma~\ref{L-supsystem}.

(iii)$\Rightarrow$(i). Let $E\subseteq X$ be a finite-dimensional operator
subsystem, set $d=\dim E$, and let $\varepsilon > 0$. Let
$\{ (x_k , x^*_k ) \}_{k=1}^d$ be an Auerbach system for $E_\sa$ 
(see for example \cite[page 335]{ER}). 
By complexifying we may view this as a biorthogonal system for $E$ with
$\| x^*_k \| \leq 2$ for each $k=1, \dots ,d$. For every
$k=1,\dots ,d$ choose a $x'_k \in X_\sa$ with $L(x'_k ) < \infty$ and
$\| x_k - x'_k \| < \varepsilon /d$. Take a $\gamma > 0$ such that for each
$k=1,\dots ,d$ we have $\gamma L( x'_k ) \leq 1$, that is, 
$\gamma x'_k \in\cE (X)$. By (iii) there is a
Lip-normed matrix algebra $(M_n , L' )$ such
that $\coedist (X, M_n ) < \gamma\varepsilon /d$. By the definition of
$\coedist$ we may view $X$ and $M_n$ as operator subsystems of an operator
system $V$ such that for each $k=1,\dots ,d$ there exists a
$y_k \in \cE (M_n )$
with $\| \gamma x'_k - y_k \| < \gamma\varepsilon /d$. Then
$$ \sum_{k=1}^d \| x_k - \gamma^{-1} y_k \| \leq \sum_{k=1}^d
\| x_k - x'_k \| +  \sum_{k=1}^d \gamma^{-1} \| \gamma x'_k - y_k \| <
2\varepsilon , $$
and thus setting $Y = \spn \{ y_1 , \dots , y_d \} \subseteq M_n$ we conclude
by Lemma~2.13.2 of \cite{Pis} that $d_\cb (E,Y) \leq
(1+4\varepsilon )(1-4\varepsilon )^{-1}$. Since $\varepsilon$ was arbitrary
this yields (i).
\end{proof}

%Denoting by $\OM_\ex$ the set of isometry classes of Lip-normed 
%$1$-exact operator systems in $\OM$, we have the following immediate 
%corollary.

%\begin{corollary}
%The set $\OM_\ex$ is closed in $\OM$.
%\end{corollary}

\begin{remark}\label{R-cb}
A perturbation argument as in the proof
of (iii)$\Rightarrow$(i) in Theorem~\ref{T-exact} shows that if a sequence
$\{ (X_n ,L_n ) \}_{n\in\Nb}$ of Lip-normed operator systems converges in
the operator Gromov-Hausdorff topology to a Lip-normed operator system
$(X,L)$, then for every finite-dimensional subspace $E\subseteq X$ there
exist, for all sufficiently large $n$, subspaces $E_n \subseteq X_n$ with 
$\dim E_n = \dim E$ such that
$\lim_{n\to\infty} d_\cb (E_n , E) = 1$. In the next section we will require
some quantitative information concerning the relationship between $\coedist$ 
and $d_\cb$ in the finite-dimensional case (see Lemma~\ref{L-dcomparison}).
\end{remark}

The above remark shows in particular that the operator space exactness 
constant (see \cite{ER,Pis}) is lower semicontinuous on $\OM$. In other words, 
if for $\lambda \geq 1$ we denote by $\OM_{\lambda\text{-}\ex}$ the set of 
isometry classes of Lip-normed $\lambda$-exact operator systems in $\OM$,
then we have:

\begin{proposition}
For each $\lambda\geq 1$, the set $\OM_{\lambda\text{-}\ex}$ is closed in 
$\OM$.
\end{proposition}

Following the notation of Section~6 of \cite{MQGHD},
for a Lip-normed operator system $(X,L)$ and $\varepsilon > 0$ we denote by
$\Afn_L (\varepsilon )$ the smallest positive integer $k$ such that there is a
Lip-normed operator system $(Y,L_Y )$ with $Y$ an operator subsystem of
the matrix algebra $M_k$ and $\sdist (X,Y) \leq\varepsilon$, and put
$\Afn_L (\varepsilon ) = \infty$ if no such
$k$ exists. By Theorems~\ref{T-equivalent} and \ref{T-exact}, for a
Lip-normed operator system $(X,L)$
we have that $X$ is $1$-exact if and only if $\Afn_L (\varepsilon )$ is
finite for all $\varepsilon > 0$. In other words,
the set $\OM_\ex$ of isometry classes of Lip-normed $1$-exact operator systems
coincides with $\mathcal{R}_\fa$,
in the notation of Section~6 in \cite{MQGHD}. We can thus restate the
compactness theorem of \cite{MQGHD} as follows.

\begin{theorem}\cite[Thm.\ 6.3]{MQGHD}\label{T-totbded}
Let $\mathcal{C}$ be a subset of $\OM_\ex$. Then $\mathcal{C}$ is
totally bounded if and only if
\begin{enumerate}
\item there is a $D>0$ such that the diameter of every element of $\mathcal{C}$
is bounded by $D$, and
\item there is a function $F : (0,\infty ) \to (0,\infty )$ such that
$\Afn_L (\varepsilon ) \leq F(\varepsilon )$ for all $(X,L)\in\mathcal{C}$.
\end{enumerate}
\end{theorem}

To establish the analogue of Theorem~\ref{T-exact} for $\nudist$
we will need the following characterization of being quasidiagonal and
exact for separable unital $C^*$-algebras. This combines
results of Voiculescu \cite{NQDH} and Dadarlat \cite{OAQDA} and the
equivalence of exactness and nuclear embeddability due to Kirchberg \cite{Kir}.

\begin{theorem}\label{T-qdexact}
A separable unital $C^*$-algebra $A$ is quasidiagonal and exact if and only if
for every finite set $\{ x_1 , \dots , x_n \}\subseteq A$ and
$\varepsilon > 0$ there is a unital $C^*$-algebra $D$, a finite-dimensional
unital $C^*$-subalgebra $B$ of $D$, elements $y_1 , \dots , y_n \in B$, and
an injective unital $^*$-homomorphism $\Phi : A \to D$ such that
$\| \Phi (x_k ) - y_k \| < \varepsilon$ for every $k=1, \dots ,n$.
\end{theorem}

\begin{theorem}\label{T-nuqdexact}
For a Lip-normed unital $C^*$-algebra $(A,L)$ the following are equivalent:
\begin{enumerate}
\item $A$ is quasidiagonal and exact,
\item for every $\varepsilon > 0$ there is a Lip-normed finite-dimensional
$C^*$-algebra $(B,L' )$ such that $\nudist (A,B) \leq\varepsilon$,
\item for every $\varepsilon > 0$ there is a Lip-normed matrix
algebra $(M_n ,L' )$ such that\linebreak $\nudist (A, M_n ) \leq\varepsilon$.
\end{enumerate}
\end{theorem}

\begin{proof}
(i)$\Rightarrow$(ii). Let $\varepsilon > 0$. Since $\cD (A)$ is
totally bounded and $\cE (A) = \cD (A) + \Rb 1$, we can find a 
finite-dimensional subspace $X$ of $A_\sa$ containing $1$ such that
$\dist_{\rm H} (\cE (A) , X\cap\cE (A)) \leq \varepsilon /3$.
Take a linear basis $\{ 1, x_1 , \dots , x_n \}$ for $X$.

By hypothesis $A$ is quasidiagonal and exact, and since it admits a Lip-norm
it must also be separable. Thus, by Theorem~\ref{T-qdexact}, given a
$\delta > 0$ we may view $A$
as a unital $C^*$-subalgebra of a unital $C^*$-algebra $D$ such that there
exists a unital finite-dimensional $C^*$-subalgebra $B$ of $D$ and a
$y_x \in B$ with $\| x - y_x \| \leq\delta$ for each
$x\in X\cap\cE (A)$. Choose a $\gamma > 0$ such that
$\gamma\rd (X) \leq \varepsilon /3$. By taking $\delta$
small enough we may assume by a standard perturbation argument that the
unital linear map $\varphi : X \to B_\sa$ defined by $\varphi (1) = 1$ and
$\varphi (x_k ) = \re\, y_{x_k}$ for $k=1, \dots ,n$,
is injective and satisfies $\| x - \varphi (x) \| \leq
\gamma \min (\| x \| , \| \varphi (x) \|)$ for all $x\in X$. Define a Lip-norm
$L''$ on $Y = \varphi (X)$ by setting $L'' (y) = L(\varphi^{-1} (y))$ for all
$y\in Y$.

Now if $x\in X\cap\cD (A)$ then setting $y = \varphi (x)$ 
we have $y\in\cE (Y)$ and $\| x - y \| \leq \gamma \| x \| \leq
\gamma \rad (X) \leq\varepsilon /3$. 
Since $\cE (A) = \cD (A) + \Rb 1$ and $\cE (Y) = \varphi (X\cap\cE (A))$, it 
follows that $\dist_{\rm H} (X\cap\cE (A) , \cE (Y)) \leq\varepsilon /3$.

By Lemma~\ref{L-supsystem} we can define
a Lip-norm $L'$ on $B$ such that $\dist_{\rm H} ( \cE (Y) ,
\cE (B)) \leq \varepsilon /3$. The triangle inequality then yields 
$\dist_{\rm H} (\cE (A) , \cE (B))
\leq \varepsilon /3 + \varepsilon /3 + \varepsilon /3 = \varepsilon$ so that
$\nudist (A,B) \leq\varepsilon$, as desired.

(ii)$\Rightarrow$(iii). Every finite-dimensional $C^*$-algebra embeds as a
unital $C^*$-subalgebra of a matrix algebra and so we can apply
Lemma~\ref{L-supsystem}.

(iii)$\Rightarrow$(i). We can apply an approximation argument similar to the
one in the proof of (iii)$\Rightarrow$(i) in Theorem~\ref{T-exact} and
appeal to Theorem~\ref{T-qdexact}.
\end{proof}

\noindent It follows from Theorems~\ref{T-exact} and \ref{T-nuqdexact} that
on the set $\CAM$ of isometry classes of Lip-normed unital $C^*$-algebras the
operator Gromov-Hausdorff topology is strictly weaker than the $\nudist$
topology.

We round out this section by giving a characterization of matrix 
approximability for the $C^*$-algebraic quantum Gromov-Hausdorff distance 
$\cqdist$ introduced in \cite{CAQGHD}. This will allow us to compare the 
$C^*$-algebraic quantum Gromov-Hausdorff and operator Gromov-Hausdorff 
topologies on $\CAM$.

The distance $\cqdist$ was defined in \cite{CAQGHD} for $C^*$-algebraic 
compact quantum metric spaces (see the beginning of Section~\ref{S-complete}),
but it only requires the Lip-norm on self-adjoint elements, and so the
definition makes sense for general Lip-normed unital $C^*$-algebras. Thus 
given Lip-normed unital $C^*$-algebras $(A,L_A )$ and $(B,L_B )$ we define
$$ \cqdist (A,B) = \inf \dist_{\rm H} \big( h^{(3)}_A (\cD (A)^{\rm m} ) ,
h^{(3)}_B (\cD (B)^{\rm m} )\big) $$
where the infimum is taken over all triples $(V, h_A , h_B )$ consisting of
a normed linear space $V$ and isometric linear maps $h_A : A\to V$ and 
$h_B : B\to V$. As in \cite{CAQGHD}, for a subset $X$ of a $C^*$-algebra $A$ 
we are using $X^{\rm m}$ to denote 
$\{ (x,y,xy) \in A\oplus A \oplus A : x,y \in X \}$ 
and for a linear map $h : V\to W$ between normed linear spaces we are 
writing $h^{(3)}$ for the induced map $V \oplus V \oplus V\to 
W \oplus W \oplus W$ between the threefold $\ell_\infty$-direct sums.

In fact given a Lip-normed unital 
$C^*$-algebra $(A,L)$ we can produce a $C^*$-algebraic compact quantum metric 
space $(A,L' )$ by setting
$$ L(a) = \sup \bigg\{ 
\frac{| \omega (a) - \sigma (a) |}{\rho_L (\omega , \sigma )} : 
\omega , \sigma\in S(A) \text{ and } \omega\neq\sigma \bigg\} $$
for all $a\in A$, in which case the set $\cD (A)$ associated to $L'$ is the 
closure of the set $\cD (A)$ associated to $L$. We can thus 
equivalently view $\CAM$ with the metric $\cqdist$ as the set of isometry 
classes of $C^*$-algebraic compact quantum metric spaces (``isometry'' has the 
same meaning in this case \cite[Def.\ 3.14]{CAQGHD}) with the metric 
$\cqdist$ as 
originally defined in \cite[Def.\ 3.3]{CAQGHD}.

Recall that a separable $C^*$-algebra $A$ is said to be an {\it MF algebra}
if it can be expressed as the inductive limit of a generalized inductive
system of finite-dimensional $C^*$-algebras \cite[Def.\ 3.2.1]{BK}.
By Theorem~3.2.2 of \cite{BK} this is equivalent to each of the following
conditions:
\begin{enumerate}
\item[(i)] there exists an injective $^*$-homomorphism 
$\Phi : A \to (\prod_{k=1}^\infty M_{n_k} )/
( \bigoplus_{k=1}^\infty M_{n_k} )$ for some sequence $\{ n_k \}$ of positive
integers,

\item[(ii)] $A$ admits an essential quasidiagonal extension by the compact 
operators $\cK$,

\item[(iii)] there exists a continuous field $(A_t )$ of $C^*$-algebras
over $\Nb \cup\{ \infty \}$ with $A_t$ finite-dimensional for every $t\in\Nb$
and $A_\infty = A$.
\end{enumerate}
If $A$ is unital then the $^*$-homomorphism $\Phi$ in (i) may be taken to
be unital, since the image of the unit under $\Phi$ can be lifted
to a projection $(p_k )$ in $\prod_{k=1}^\infty M_{n_k}$, yielding
an injective unital $^*$-homomorphism from $A$ to 
$(\prod_{k=1}^\infty p_k M_{n_k} p_k )/
( \bigoplus_{k=1}^\infty p_k M_{n_k} p_k )$. Thus, in view of the proof
of \cite[Prop.\ 2.2.3]{BK}, we may also in this case take 
the fibre $C^*$-algebras to be unital and the unit section to be continuous 
in the continuous field in (iii). In the proof of Theorem~\ref{T-cqmatrix} 
below we will implicitly use this unital version of (iii) as a 
characterization of being an MF algebra for separable unital $C^*$-algebras.

Note that a Lip-normed unital $C^*$-algebra is automatically
separable as a $C^*$-algebra.

\begin{lemma}\label{L-ctsfield}
Let $( A_t )$ be a continuous field of $C^*$-algebras 
over $\Nb \cup \{ \infty \}$ with separable unital fibres such that the unit 
section is continous. Let $L$ be a Lip-norm on $A_\infty$. Then 
for each $t\in\Nb$ there is a Lip-norm $L_t$ on $A_t$ so that $(A_t , L_t )$
forms a continuous field of Lip-normed unital $C^*$-algebras
over $\Nb \cup \{ \infty \}$ (with the same meaning as that given in
\cite[Def.\ 5.1]{CAQGHD}) and
$\lim_{t\to\infty} \GHdist (\cD (A_t ) , \cD (A_\infty )) = 0$.
\end{lemma}

\begin{proof}
By \cite[Cor.\ 2.8]{GGh} we can find a normed linear space $V$
containing each $A_t$ isometrically such that for every 
continuous section $f$ the map $t\mapsto f_t$ is continuous at $t = \infty$.
Let $\varepsilon > 0$. Take a finite-dimensional subspace $X$ of $A_\sa$ such 
that $\Hdist (\cD (A_\infty ) , X\cap \cD (A_\infty )) \leq\varepsilon /3$. 
As in the proof of (i)$\Rightarrow$(ii) in Theorem~5.6, by a standard 
perturbation argument we can construct, for any sufficiently large $t\in\Nb$, 
a unital linear map $\varphi : X \to (A_t )_\sa$ such that by taking
the Lip-norm $L''$ on $Y = \varphi (X)$ defined by $L'' (y) = 
L(\varphi^{-1} (y))$ we get $\Hdist (X \cap \cD (A), \cD (Y)) \leq
\varepsilon /3$. By Lemma~\ref{L-supsystem} there is a Lip-norm
$L'$ on $A_t$ such that $\Hdist (\cD (Y), \cD (B)) \leq \varepsilon /3$.
The triangle inequality then yields 
$$ \GHdist (\cD (A_\infty ) , \cD (A_t )) \leq 
\Hdist (\cD (A_\infty ) , \cD (A_t )) \leq \varepsilon /3 + \varepsilon /3
+ \varepsilon /3 = \varepsilon . $$
We can thus apply this procedure to produce a Lip-norm $L_t$ on $A_t$ for
each $t\in\Nb$ so that 
$\lim_{t\to\infty} \GHdist (\cD (A_t ) , \cD (A_\infty )) = 0$, and it is
readily checked that $(A_t , L_t )$ forms a continuous field of
Lip-normed unital $C^*$-algebras.
\end{proof}

\begin{theorem}\label{T-cqmatrix} 
Let $(A,L)$ be a Lip-normed unital $C^*$-algebra. Then the 
following are equivalent:
\begin{enumerate}
\item $A$ is an MF algebra,
\item for every $\varepsilon > 0$ there is a finite-dimensional
Lip-normed unital $C^*$-algebra $(B,L' )$ such that 
$\cqdist (A,B) \leq\varepsilon$,
\item for every $\varepsilon > 0$ there is a Lip-normed matrix algebra 
$(M_n ,L' )$ such that\linebreak $\cqdist (A, M_n ) \leq\varepsilon$.
\end{enumerate}
\end{theorem}

\begin{proof}
The implication (i)$\Rightarrow$(ii) is a consequence of 
Lemma~\ref{L-ctsfield} above and Theorem~5.2 of \cite{CAQGHD}, while 
(ii)$\Rightarrow$(iii) follows from Lemma~\ref{L-supsystem} and 
(iii)$\Rightarrow$(i) from Proposition~5.4 of \cite{CAQGHD}.
\end{proof}

\noindent Theorem~\ref{T-cqmatrix} shows that on the set $\CAM$ of isometry 
classes of Lip-normed unital $C^*$-algebras the $\cqdist$ topology is strictly
weaker than the $\nudist$ topology and is neither weaker nor stronger than
the operator Gromov-Hausdorff topology.

Finally, we remark that the compactness theorem recorded above as 
Theorems~\ref{T-totbded} also holds, with the appropriate substitutions, for 
$\nudist$ and $\cqdist$. In these cases we take $\Afn_L (\varepsilon )$ 
for a Lip-normed unital $C^*$-algebra $(A,L)$ to be the
smallest positive integer $k$ such that there is a Lip-normed
matrix algebra $(M_k ,L' )$ with $\nudist (A,M_k ) \leq\varepsilon$
(resp.\ $\cqdist (A,M_k ) \leq\varepsilon$) and setting 
$\Afn_L (\varepsilon ) = \infty$ if no such $k$ exists. Then 
for $\nudist$ we should replace $\OM_\ex$ by the set of isometry classes of
Lip-normed unital quasidiagonal exact $C^*$-algebras, as 
Theorem~\ref{T-nuqdexact} shows, while for $\cqdist$ we should instead
substitute the set of isometry classes of Lip-normed unital MF $C^*$-algebras, 
by Theorem~\ref{T-cqmatrix}. To establish the compactness theorems, we use the
observation that the set of isometry classes 
of Lip-normed matrix algebras $(M_k , L)$ for a fixed $k$ and with a 
fixed upper bound on the radius is totally bounded under both $\nudist$ and 
$\cqdist$. For $\nudist$ this follows from the fact that, given Lip-norms
$L_1$ and $L_2$ on a matrix algebra $M_k$, the quantity
$\nudist ((M_k , L_1 ) , (M_k , L_2 ))$ is bounded above by  
$\Hdist (\cE (M_k , L_1 ) , \cE (M_k , L_2 ))$, which coincides with the 
Hausdorff distance in $M_k / \Cb 1$ between the images
of $\cE (M_k , L_1 )$ and $\cE (M_k , L_2 )$ under the quotient. 
Compare the proof of Theorem~6.3 in \cite{MQGHD}.

\section{Nonseparability}\label{S-nonseparable}

We write $\OS_n$ for the set of $n$-dimensional operator spaces, with two
such operator spaces being considered the same if they are completely
isometric. The subset of $n$-dimensional Hilbertian operator
spaces (i.e., operator spaces which are isometric to
$\ell^n_2$ as normed spaces) will be denoted $\HOS_n$.
We write $\OM_n$ for the set of isometry classes of Lip-normed
$n$-dimensional operator systems equipped with the operator Gromov-Hausdorff
topology, i.e., the $\coedist$ (or equivalently $\sdist$) topology. Note
that for $n\geq 2$ the set $\OM_n$ is not closed in $\OM$ since it is possible 
to have dimension collapse (as happens already for finite metric spaces). 
However, $\bigcup_{1\leq n \leq m} \OM_n$ is closed in $\OM$ for each 
$m\in\Nb$, as can be gathered from Remark~\ref{R-cb}. 

Let $(X,L)$ be a Lip-normed operator system. We denote by $\tilde{L}$ the norm
on the Banach space quotient $X_\sa / \Rb 1$ induced from $L$. Suppose that
$X$ is of some finite dimension strictly greater than one.
Then the formal identity map $I_X : (X_\sa / \Rb 1, \tilde{L} ) \to
(X_\sa / \Rb 1, \| \cdot \| )$ is an isomorphism, and we can define
$\mu_X = \| I_X^{-1} \|$. 
%$$ \mu_X = \max (\| I_X \|^{-1} , \| I_X^{-1} \| ) . $$
%The quantity $\| I_X \|$ is the radius of $(X,L)$ \cite[Prop.\ 2.2]{MSS},
%while it is easily checked that
%$$ \| I_X^{-1} \| = \sup_{x\in X_\sa \setminus \{ 0 \}}
%\frac{L(x)}{\| x \|} . $$
In the case that $X = C(F)$ for a finite metric space $F$,
$\mu_X$ is equal to the inverse of half of the smallest distance
between any two points of $F$.

\begin{lemma}\label{L-dcomparison}
Let $(X,L_X )$ and $(Y,L_Y )$ be Lip-normed operator systems of some finite
dimension $n\geq 2$. Set $\kappa_{X,Y} = 2n \min (\mu_X , \mu_Y )
\hspace*{0.2mm}\coedist (X,Y)$.
If $0 \leq \kappa_{X,Y} < 1$ then we have
$$ d_\cb (X,Y) \leq \frac{1 + \kappa_{X,Y}}{1 - \kappa_{X,Y}} . $$
%$$ \frac{1 + d_\cb (X,Y)}{1 - d_\cb (X,Y)} \leq n \min (\mu (X), \mu (Y))
%\coedist (X,Y) . $$
\end{lemma}

\begin{proof}
Without loss of generality we may assume that $\mu_X \leq \mu_Y$.
Let $\{ (x_k , x^*_k ) \}_{k=1}^n$ be an Auerbach system for $X_\sa$ (see for 
example \cite[page 335]{ER}). By
complexifying we regard this as a biorthogonal system for $X$ with
$\| x^*_k \| \leq 2$ for each $k=1, \dots ,n$. For each $k=1, \dots ,n$ we
have $L(x_k ) \leq \mu_X \| x_k \| = \mu_X$, that is, 
$\mu_X^{-1} x_k \in\cE (X)$. By the definition of
$\coedist$ we may view $X$ and $Y$ as operator subsystems of an
operator system $Z$ such that for each $k=1, \dots ,n$ there is a
$y_k \in\cE (Y)$ with $\| \mu_X^{-1} x_k - y_k \| \leq\coedist (X,Y)$.
We then have
$$ \sum_{k=1}^n \| x^*_k \| \| x_k - \mu_X y_k \| \leq
2 \sum_{k=1}^n \mu_X \| \mu_X^{-1} x_k - y_k \| \leq 2n \mu_X
\coedist (X,Y) , $$
and since $X$ and $Y$ are of the same dimension we conclude by Lemma~2.13.2 of
\cite{Pis} that $d_\cb (X,Y) \leq (1 + \kappa_{X,Y} )(1 - \kappa_{X,Y} )^{-1}$,
as desired.
\end{proof}

\begin{theorem}\label{T-nonseparable}
For every $n\geq 7$ there is a nonseparable (and in particular
non-totally-bounded) subset of $\OM_n$ whose elements are all isometric to
each other as compact quantum metric spaces.
\end{theorem}

\begin{proof}
For each $V\in\HOS_3$ we take a representation of $V$ as an operator
subspace of a unital $C^*$-algebra $A$ and define the $7$-dimensional
operator system
$$ X_V = \left\{ \left[
\begin{matrix} \lambda 1 & a \\ b^* & \lambda 1 \end{matrix} \right]
\in M_2 (A) : \lambda\in\Cb\text{ and } a,b\in V \right\} . $$
Note that we have a completely isometric
embedding $V\hookrightarrow X_V$ given by
$a\mapsto \big[ \begin{smallmatrix} 0 & a \\ 0 & 0 \end{smallmatrix} \big]$.
The self-adjoint elements of $X_V$ are those of the form
$\big[ \begin{smallmatrix} \lambda 1 & a \\ a^* & \lambda 1 \end{smallmatrix}
\big]$ for $\lambda\in\Rb$ and $a\in V$, and by Lemma~3.1 of \cite{Pau}
such an element is positive if and only if $\lambda\geq \| a \|$. Thus if
we define $Y_V$ as the operator system direct sum of $X_V$ and the
commutative $C^*$-algebra $\Cb^{n-7}$ we see that the
operator systems $Y_V$ for $V\in\HOS_3$ are all order-isomorphic to each
other. Moreover if we define a Lip-norm $L_V$ on $(Y_V )_\sa$ by taking
$L_V (y)$ to be the norm of the image of $y$ under the quotient map
$(Y_V )_\sa \to (Y_V )_\sa / \Rb 1$ then the
Lip-normed operator systems $(Y_V , L_V )$ for
$V\in\HOS_3$ are all isometric to each other as compact quantum metric
spaces. Denote by $\Theta$ the subset of $\OM_n$ consisting of the
$(Y_V , L_V )$ and by $\Gamma$ the subset of $\OS_n$ consisting of the
$Y_V$.

We claim that $\Gamma$ is nonseparable. Suppose that this is not the case.
Then since for any integers $s\geq r\geq 1$ the set of $r$-dimensional
subspaces of a given $s$-dimensional operator space is compact in $\OS_r$
(this can be shown using Lemma~2.13.2 of \cite{Pis}) we infer the separability
of the subset of $\OS_3$ consisting of all $3$-dimensional operator spaces
which appear as a subspace of some operator system in $\Gamma$.
But this subset of $\OS_3$ contains $\HOS_3$, which is
nonseparable \cite[Remark 2.4]{JP}, producing a contradiction.
Thus $\Gamma$ is nonseparable, and so by Lemma~\ref{L-dcomparison}
we conclude that $\Theta$ is nonseparable, as desired.
\end{proof}

\noindent Theorem~\ref{T-nonseparable} shows that, in contrast 
to the order-unit case \cite[Thm.\ 13.5]{GHDQMS}\cite[Thm.\ 5.5]{OUQGHD}, 
there can be no compactness theorem for $\coedist$ or $\sdist$ which at the
``topological'' level makes 
reference only to the state space or norm structure 
(cf.\ Theorem~\ref{T-totbded}). It also follows that for every 
$n\geq 7$ and $D>0$ the set of isometry classes of Lip-normed $n$-dimensional 
operator systems of diameter at most $D$ is not separable and, in particular,
not totally bounded, thus answering Question~6.5 of \cite{MQGHD} for $n\geq 7$.

\section{Generic complete order structure}\label{S-generic}

In Theorem~\ref{T-OMexgeneric} we describe the type of complete order structure
possessed by a generic element in $\OM_\ex$ under operator Gromov-Hausdorff distance,
where $\OM_\ex$ is viewed as consisting of Lip-normed $1$-exact operator systems $(X,L)$ 
with $X$ complete and $L$ closed. Note that $\OM_\ex$ is a separable closed subset of the
complete metric space $\OM$ by Theorem~\ref{T-exact}, and so $\OM_\ex$ is a Baire space.
%As a point of comparison, we also prove an analogous result for Lip-normed quasidiagonal 
%exact unital $C^*$-algebras under $\nudist$, with the conclusion that 
%the generic $C^*$-algebraic structure is that of the universal UHF algebra 
%(Theorem~\ref{T-CMqdegeneric}). 

%We will view elements of $\OM_\ex$ as Lip-normed operator systems $(X,L)$ with
%$X$ complete and $L$ closed, 
%identifying two such Lip-normed operator systems if 
%they are isometric (in which case the operator systems are unitally completely 
%order isomorphic).

\begin{lemma}\label{L-ucbucoe}
Let $d\in\Nb$ and $\varepsilon > 0$. Then there exists a $\delta > 0$ such
that whenever $n\in\Nb$ and $\varphi : M_d \to M_n$ is an injective unital
linear map with $\max (\| \varphi \|_\cb , \| \varphi^{-1} \|_\cb ) < 
1 + \delta$, there exists a unital complete order embedding
$\psi : M_d \to M_n$ with $\| \psi - \varphi \| < \varepsilon$.
\end{lemma}

\begin{proof}
Suppose that the lemma is not true. Then there is an 
$\varepsilon > 0$ such that for each
$k\in\Nb$ there exists an injective unital linear map $\varphi_k$ from
$M_d$ to a matrix algebra $M_{n_k}$ such that
$\max (\| \varphi_k \|_\cb , \| \varphi_k^{-1} \|_\cb ) < 1 + 1/k$ and
$\| \psi - \varphi \| \geq \varepsilon$ for every unital complete order
embedding $\psi : M_d \to M_{n_k}$. Define the map $\varphi : M_d \to
(\prod_{k=1}^\infty M_{n_k} )/(\bigoplus_{k=1}^\infty M_{n_k} )$ by 
$\varphi (x) = \pi ((\varphi_k (x))_k )$ where $\pi : \prod_k M_{n_k} \to
(\prod_{k=1}^\infty M_{n_k} )/(\bigoplus_{k=1}^\infty M_{n_k} )$ is the 
quotient map. Then $\varphi$ is a unital complete isometry and hence
a complete order embedding \cite[Cor.\ 5.1.2]{ER}.

Let $\{ e_{ij} \}_{1\leq i,j \leq d}$ be the set of standard matrix units 
for $M_d$. By Proposition~4.2.8 of \cite{BK} (or rather the unital version
which follows from the same proof) there is a 
$\delta > 0$ such that whenever $\gamma$ is a u.c.p.\ map from $M_d$ to
a finite-dimensional $C^*$-algebra $B$ with 
$\| \gamma (e_{12} )\gamma (e_{23} ) \cdots \gamma (e_{d-1,d} ) \| \geq 
1-\delta$, there is a unital complete order embedding $\gamma' : M_d \to B$
with $\| \gamma' - \gamma \| < \varepsilon /2$. 

By the Choi-Effros lifting theorem \cite{CPLP} there is a u.c.p.\ map 
$\theta : M_d \to\prod_{k=1}^\infty M_{n_k}$ such that $\pi\circ\theta = \varphi$. 
By the compactness of the unit ball of $M_d$ we can find a $j\in\Nb$ such that 
$\| \pi_j \circ\theta - \varphi_j \| < \varepsilon /2$, where
$\pi_j : \prod_{k=1}^\infty M_{n_k} \to M_{n_j}$ is
the projection map onto the $j$th coordinate. 
Since $\theta (e_{12} ) \theta (e_{23} ) \cdots \theta (e_{d-1,d} )$ 
is a lift of $\varphi (e_{12} )\varphi (e_{23} ) \cdots \varphi (e_{d-1,d} )$
under $\pi$ and 
\[ \| \varphi (e_{12} )\varphi (e_{23} ) \cdots \varphi (e_{d-1,d} ) \| \geq 
\| \varphi (e_{12} e_{23} \cdots e_{d-1,d} ) \| = \| \varphi (e_{1d} ) \| = 1 \]
by the supermultiplicativity of $\varphi$ \cite[Prop.\ 4.2.5]{BK}, 
we may assume that $j$ is large enough so that we additionally have 
\[ \| \pi_j (\theta (e_{12} ))\pi_j (\theta (e_{23} )) \cdots 
\pi_j (\theta (e_{d-1,d} )) \| = 
\| \pi_j ( \theta (e_{12} ) \theta (e_{23} ) \cdots \theta (e_{d-1,d} )) \|
\geq 1 - \delta . \]
Because $\pi_j \circ\theta$ is a u.c.p.\ map, it follows by our choice of $\delta$
that there is a unital complete order embedding $\psi : M_d \to M_{n_j}$  
with $\| \psi - \pi_j \circ\theta \| < \varepsilon /2$. But then
\[ \| \psi - \varphi_j \| \leq \| \psi - \pi_j \circ\theta \| 
+ \| \pi_j \circ\theta - \varphi_j \| < \varepsilon /2 + \varepsilon /2 = 
\varepsilon , \]  
producing a contradiction.
\end{proof}

\begin{lemma}\label{L-perturb}
Let $(M_d ,L)$ be a Lip-normed matrix algebra and let
$\varepsilon > 0$. Then there is a $\delta > 0$ such that whenever
$(X ,L_X )$ is a Lip-normed operator system, $A$ is a unital 
$C^*$-algebra, and $\beta : M_d \to A$ and $\gamma : X \to A$ are unital 
complete order embeddings for which 
$\Hdist (\beta (\cE (M_d )), \gamma (\cE (X))) < \delta$, there exists an 
injective unital linear map $\varphi : M_d \to X$ such that 
$\max (\| \varphi \|_\cb , \| \varphi^{-1} \|_\cb ) \leq 1 + \varepsilon$ and
$\| \beta - \gamma\circ\varphi \| < \varepsilon$.
\end{lemma}

\begin{proof}
We may assume that $d>1$. Letting $\{ e_{ij} \}_{1\leq i,j\leq d}$ be the set 
of standard matrix units for $M_d$, for $i,j=1,\dots ,d$ we set 
$a_{ij} = e_{ij} + e_{ji}$ if $i>j$, $a_{ij} = i(e_{ij} - e_{ji} )$ if 
$i<j$, and $a_{ij} = e_{ii}$ if $i=j$. Then $\{ a_{ij} \}_{1\leq i,j\leq d}$ is 
an Auerbach basis for $M_d$. Since $M_d$ is finite-dimensional 
there exists a $K>0$ such that
$L(a) \leq K \| a \|$ for all self-adjoint $a\in M_d$.
Let $0 < \eta < \min (\varepsilon ,1)$, and suppose that 
$(X ,L_X )$ is a Lip-normed finite-dimensional $C^*$-algebra, $A$ is a unital 
$C^*$-algebra, and $\beta : M_d \to A$ and $\gamma : X \to A$ are unital 
complete order embeddings for which 
\[ \Hdist (\beta (\cE (M_d )), \gamma (\cE (X)) < 
\frac{\eta}{Kd^3} . \] 
Then the image under $\beta$ of the norm unit ball of $(M_d )_\sa$ is contained
in $K\beta (\cE (M_d ))$, and so for all $i,j=1,\dots ,d$
we can find an $x_{ij} \in X_\sa$ such that 
$\| \beta (a_{ij} ) - \gamma (x_{ij} ) \| < d^{-3} \eta$. Redefining 
$x_{11}$ as $1-\sum_{i=2}^d x_{ii}$, we then have
$\| \beta (a_{ij} ) - \gamma (x_{ij} ) \| < d^{-2} \eta$ for all $1\leq i,j\leq d$,
and the unital linear map $\varphi : M_d \to B$ 
determined by $\psi (a_{ij} ) = x_{ij}$ satisfies 
$\| \beta - \gamma\circ\psi \| < \eta < \varepsilon$. Moreover by Lemma~2.13.2 
of \cite{Pis} $\psi$ is injective and $\| \psi \|_\cb \leq 1+\eta$ and 
$\| \psi^{-1} \|_\cb \leq (1-\eta )^{-1}$. We can thus find an $\eta$ small enough 
as a function of $\varepsilon$ to ensure that $\max (\| \varphi \|_\cb ,
\| \varphi^{-1} \|_\cb ) \leq 1 + \varepsilon$ and then take 
$\delta = \eta (Kd^3 )^{-1}$ to obtain the lemma.
\end{proof}

\begin{lemma}\label{L-perturbmatrix}
Let $(M_d ,L)$ be a Lip-normed matrix algebra and let
$\varepsilon > 0$. Then there is a $\delta > 0$ such that whenever
$(M_n ,L' )$ is a Lip-normed matrix algebra, $A$ is a unital 
$C^*$-algebra, and $\beta : M_d \to A$ and $\gamma : M_n \to A$ are unital 
complete order embeddings for which 
$\Hdist (\beta (\cE (M_d )), \gamma (\cE (M_n ))) < \delta$, there exists a 
unital complete order embedding $\varphi : M_d \to M_n$ satisfying
$\| \beta - \gamma\circ\varphi \| < \varepsilon$.
\end{lemma}

%%%%%%%%%%%%%%%%%%%%%%
% check E/D
%%%%%%%%%%%%%%%%%%%%%%

\begin{proof}
By Lemmas~\ref{L-ucbucoe} and \ref{L-perturb}, if $\delta$ is sufficiently small 
as a function of $d$ and $\varepsilon$ then whenever $(M_n ,L' )$ is a Lip-normed 
matrix algebra, $A$ is a unital $C^*$-algebra, and 
$\beta : M_d \to A$ and $\gamma : M_n \to A$ are unital complete order embeddings for 
which\linebreak $\Hdist (\beta (\cE (M_d )), \gamma (\cE (M_n ))) < \delta$,
there exist an injective unital linear map 
$\psi : M_d \to M_n$ with $\| \beta - \gamma\circ\psi \| < \varepsilon /2$
and a unital complete order embedding
$\varphi : M_d \to M_n$ with $\| \psi - \varphi \| < \varepsilon /2$,
in which case 
\[ \| \beta - \gamma\circ\varphi \| \leq \| \beta - \gamma\circ\psi \|
+ \| \gamma \| \| \psi - \varphi \| < \varepsilon . \]
\end{proof}

\begin{lemma}\label{L-matrixdense}
Let $I$ be an infinite subset of $\Nb$, and let $\Lambda_I$ be the subset
of $\OM_\ex$ consisting of all Lip-normed matrix algebras $(M_d ,L)$ such 
that $d\in I$. Then $\Lambda_I$ is dense in $\OM_\ex$.
\end{lemma}

\begin{proof}
For positive integers $d\leq n$ there always exists a unital complete order
embedding $\varphi : M_d \to M_n$. For example, take a state $\sigma$ on
$M_d$, a rank $d$ projection $p\in M_n$, and a $^*$-isomorphism
$\Phi : M_d \to pM_n p$, and define $\psi (x) = \Phi (x) + \sigma (x) (1-p)$
for all $x\in M_d$. Using this fact in conjunction with 
Lemma~\ref{L-supsystem} and Theorem~\ref{T-exact} yields the lemma.
\end{proof}

\begin{theorem}\label{T-OMexgeneric}
Let $I$ be an infinite subset of $\Nb$. Then, with respect to the operator
Gromov-Hausdorff topology, a generic element of $\OM_\ex$ is, as an operator
system, unitally completely order isomorphic to an operator system inductive 
limit $\lim\limits_{\longrightarrow} (M_{n_k} ,\varphi_k )$ over $\Nb$ 
where $\{ n_k \}_k$ is a strictly increasing sequence in $I$ and each 
$\varphi_k$ is a unital complete order embedding.
\end{theorem}

\begin{proof}
For every Lip-normed matrix algebra $(M_d ,L)$ and every 
$\varepsilon > 0$ we take a $\delta (d,L,\varepsilon ) \in (0,\varepsilon )$ 
that works for both Lemma~\ref{L-perturb} and Lemma~\ref{L-perturbmatrix}
and define $\Gamma_\varepsilon (M_d , L)$ to be the set
of all $(X,L_X )$ in $\OM_\ex$ such that
\[ \coedist ((M_d ,L) , (X,L_X )) < \delta (d,L,\varepsilon )/2 . \]
For every subset $I\subseteq\Nb$ and $\varepsilon > 0$ we define the open
subset $\Theta (I,\varepsilon )$ of $\OM_\ex$ as the union of the sets
$\Gamma_\varepsilon (M_d , L)$ over all Lip-normed matrix algebras
$(M_d ,L)$ with $d\in I$.
By Lemma~\ref{L-matrixdense}, for all infinite subsets $I\subseteq\Nb$ and
$\varepsilon > 0$ the set $\Theta (I,\varepsilon )$ is dense in $\OM_\ex$, 
and so the countable intersection
\[ R := \bigcap_{\substack{J\subseteq I \\ I\setminus J \text{ finite}}}
\bigcap_{k\in\Nb} \Theta (J, 1/k) \]
is a dense $G_\delta$ subset of $\OM_\ex$. Letting $(X,L)$ be an element of
$R$, it thus suffices to show that $X$ can be expressed as an inductive limit
of the type described in the theorem statement. 

Let $\{ \varepsilon_k \}_k$ be a summable sequence of positive
real numbers. By the definition of $R$ there is a sequence 
$\{ (M_{n_k} , L_k ) \}_k$ of Lip-normed matrix algebras such that 
$\{ n_k \}_k$ is a strictly increasing sequence in $I$ and
\[ \coedist ((M_{n_k} ,L_k ) , (X,L)) < \delta (n_k , L, \varepsilon_k )/2 . \]
By Corollary~\ref{C-embedamalgsumosy} and Lemma~\ref{L-directlimitosy} there is 
a unital $C^*$-algebra $A$ containing $X$ as an operator subsystem and unital 
complete order embeddings $\beta_k : M_{n_k} \to A$ such that
$\Hdist (\beta_k (\cE (M_{n_k} )), \cE (X)) < \delta (n_k , L, \varepsilon_k )/2$
for all $k\in\Nb$. By passing to a subsequence and relabeling if necessary we may 
assume that $\delta (n_k , L, \varepsilon_k )$ decreases with $k$, and so by
the triangle inequality for Hausdorff distance and the definition of 
$\delta (n_k , L, \varepsilon_k )$ there exist for each $k\in\Nb$ 
a unital complete order embedding $\varphi_k : M_{n_k} \to M_{n_{k+1}}$ and an 
injective unital linear map $\theta_k : M_{n_k} \to X$ such that
$\| \beta_k - \beta_{k+1} \circ\varphi_k \| < \varepsilon_k$,
$\max (\| \theta_k \|_\cb , \| \theta_k^{-1} \|_\cb ) < 1 + \varepsilon_k$, and 
$\| \beta_k - \theta_k \| < \varepsilon_k$. Using the summability
of $\{ \varepsilon_k \}_k$, a simple estimate shows that for every
$k\in\Nb$ and $x\in M_{n_k}$ the sequence 
$\{ (\theta_{k+j} \circ\varphi_{k+j-1} \circ\cdots\circ\varphi_{k+1} \circ
\varphi_k )(x) \}_j$ in $X$ is Cauchy; denote by $\psi_k (x)$ its limit.
Since $\max (\| \theta_k \|_\cb , \| \theta_k^{-1} \|_\cb )\to 1$ as $k\to\infty$,
each of the resulting unital linear maps $\psi_k : M_{n_k} \to X$ is completely isometric 
(cf.\ the proof of Lemma~2.13.2 in \cite{Pis}).
These maps are compatible with the inductive system $\{ (M_{n_k} ,\varphi_k ) \}_k$ 
and thus give rise to a unital map 
$\psi : \lim\limits_{\longrightarrow} (M_{n_k} ,\varphi_k ) \to X$ which 
is surjective and completely isometric and hence a unital complete order 
isomorphism \cite[Cor.\ 5.1.2]{ER}, completing the proof.
\end{proof}

\begin{remark}
An operator system which can be expressed as an inductive limit as in 
Theorem~\ref{T-OMexgeneric} for a given $I$ is far from being
unique. Indeed it can be seen
from the results and arguments of \cite{BK2} (see Proposition~5.12 and
Theorem~5.13 therein) that the unital $C^*$-algebras that as operator systems
can be so expressed for a given $I$ are precisely the infinite-dimensional
unital prime strong NF algebras, and $C^*$-algebras are determined up to
$^*$-isomorphism by their complete order structure.
\end{remark}

\end{document}